\documentclass[12pt]{amsart}

\usepackage{hyperref}
\usepackage{amsmath, amsthm, amsfonts, amssymb}
\usepackage{mathrsfs}
\usepackage{mathtools}
\usepackage{upgreek}
\usepackage{listings}

\makeatletter
\@namedef{subjclassname@2020}{%
  \textup{2020} Mathematics Subject Classification}
\makeatother

\theoremstyle{plain}
\newtheorem{thm}{Theorem}[section]
\newtheorem*{thm*}{Theorem}{}
\newtheorem{coro}[thm]{Corollary}
\newtheorem{prop}[thm]{Proposition}

\newtheorem{lemm}[thm]{Lemma}

\theoremstyle{definition}
\newtheorem{deff}[thm]{Definition}
\newtheorem{examp}[thm]{Example}

\theoremstyle{remark}
\newtheorem{rema}[thm]{Remark}

\newtheorem{conj}[thm]{Conjecture}

\newcommand\twoscript[2]{\substack{{#1} \\ {#2}}}

\newcommand\rmi{\mathrm{i}}
\newcommand\rme{\mathrm{e}}

\newcommand\legendre[2]{\genfrac{(}{)}{}{}{#1}{#2}}
\newcommand\tbtmat[4]{\left(\begin{smallmatrix}{#1} & {#2} \\ {#3} & {#4}\end{smallmatrix}\right)}
\newcommand\tbtMat[4]{\begin{pmatrix}{#1} & {#2} \\ {#3} & {#4}\end{pmatrix}}
\newcommand*\abs[1]{\lvert#1\rvert}
\newcommand\etp[1]{\mathfrak{e}\left(#1\right)}

\newcommand\diff{\,\mathrm{d}}

\newcommand\ord{\mathop{\mathrm{ord}}}

\newcommand\numZ{\mathbb{Z}}
\newcommand\numQ{\mathbb{Q}}

\newcommand\numC{\mathbb{C}}
\newcommand\projQ{\mathbb{P}^1(\mathbb{Q})}
\newcommand\numgeq[2]{\mathbb{#1}_{\geq #2}}

\newcommand\slZ{\mathrm{SL}_2(\mathbb{Z})}

\newcommand\slR{\mathrm{SL}_2(\mathbb{R})}

\newcommand\pslZ{\mathrm{PSL}_2(\mathbb{Z})}

\newcommand\uhp{\mathfrak{H}}

\newcommand\elesltRaD[6]{\left(\left(\begin{smallmatrix}{#1} & {#2} \\ {#3} & {#4}\end{smallmatrix}\right),{#5}\left({#3} \tau+{#4} \right)^{\frac{1}{#6}}\right)}
\newcommand\elesltRaDs[6]{\left(\left(\begin{smallmatrix}{#1} & {#2} \\ {#3} & {#4}\end{smallmatrix}\right),{#5}\right)}
\newcommand\Dcover[2]{\widetilde{#1^{#2}}}

\begin{document}


\baselineskip=17pt


\title[Eisenstein series of rational weights]{Holomorphic Eisenstein series of rational weights and special values of Gamma function}

\author{Xiao-Jie Zhu}
\address{School of Mathematical Sciences\\ Tongji University\\
1239 Siping Road\\ Shanghai, P.R. China}
\email{zhuxiaojiemath@outlook.com}
\urladdr{https://orcid.org/0000-0002-6733-0755}

\date{}

\begin{abstract}
We give all possible holomorphic Eisenstein series on $\Gamma_0(p)$, of rational weights greater than $2$, and with multiplier systems the same as certain rational-weight eta-quotients at all cusps. We prove they are modular forms and give their Fourier expansions. We establish four sorts of identities that equate such series to rational-weight eta-quotients. As an application, we give series expressions of special values of Euler Gamma function at any rational arguments. These expressions involve exponential sums of Dedekind sums.
\end{abstract}



\maketitle

\tableofcontents

\section{Introduction}
\label{sec:Introduction}

We present all possible holomorphic Eisenstein series (in one complex variable $\tau$) on $\Gamma_0(p)$, of rational weights greater than $2$, and with multiplier systems the same as certain rational-weight eta-quotients. These functions involve exponentials of Dedekind sums, which are defined by
\begin{equation*}
s(h,k)=\sum_{r \bmod k}\Big(\Big(\frac{r}{k}\Big)\Big)\cdot\Big(\Big(\frac{hr}{k}\Big)\Big),\quad h\in\numZ,\,k\in\numgeq{Z}{1}.
\end{equation*}
In the above definition, $((x))=x-[x]-1/2$ if $x$ is a non-integral real number, and $((x))=0$ if $x$ is an integer. Let $p$ be a prime, $r_1$ and $r_p$ be rational numbers such that $(r_1+pr_p)/24$ is an integer. Suppose $k \in \frac{r_1+r_p}{2}+2\numZ$ is greater than $2$. Our first main theorem says that, the following Eisenstein series
\begin{multline}
\label{eq:Eisinfty}
E_{k,p}^{\rmi\infty}(\tau;r_1,r_p)=1+\sum_{\twoscript{p \mid c > 0,d\in\numZ}{\gcd(c,d)=1}}(c\tau+d)^{-k} \\
\times\etp{-r_1\left(\frac{a+d}{24c}+\frac{1}{2}s(-d,c)-\frac{1}{8}\right)-r_p\left(\frac{a+d}{24c/p}+\frac{1}{2}s(-d,c/p)-\frac{1}{8}\right)},
\end{multline}
where for each pair $(c,d)$, $a$ is any integer satisfying $ad \equiv 1 \pmod c$, is a holomorphic modular form of weight $k$, on the group $\Gamma_0(p)$ and with a multiplier system the same as the eta-quotient $\eta^{r_1}(\tau)\eta^{r_p}(p\tau)$.

We collect some common notations here. The symbol $\etp{x}$ denotes $\exp(2\uppi\rmi x)$. The variable $\tau$ is always assumed to take values from the upper half plane $\uhp = \{z \in \numC\colon \Im(z) > 0\}$ and $q^n$ denotes $\etp{n\tau}$. The symbol $\Gamma_0(N)$ denotes the congruence subgroup of $\slZ$ consisting of matrices $\tbtmat{a}{b}{c}{d}$ with $N \mid c$, and $\eta(\tau)$ denotes the Dedekind eta function defined by
\begin{equation*}
\eta(\tau)=q^{1/24}\prod_{n \in \numgeq{Z}{1}}(1-q^n),\qquad\tau\in\uhp.
\end{equation*}
For the meaning of rational-weight modular forms, of rational powers of Dedekind eta function, and of Eisenstein series of rational weights the reader may scan through Section \ref{sec:Modular forms of rational weights}, \ref{sec:Eisenstein series of rational weights} and \ref{sec:Eta quotients of rational weights}. 

Besides $E_{k,p}^{\rmi\infty}(\tau;r_1,r_p)$, we state and prove a parallel conclusion on the Eisenstein series $E_{k,p}^{1}(\tau;r_1,r_p)$ at another cusp of $\Gamma_0(p)$, but this time it is required that $(pr_1+r_p)/24$ is an integer, not $(r_1+pr_p)/24$. For details, see Theorem \ref{thm:mainEis}.

The second result concerns Fourier expansions of $E_{k,p}^{\rmi\infty}(\tau;r_1,r_p)$ and $E_{k,p}^{1}(\tau;r_1,r_p)$. The strategy of the calculation is not new: it is typical in analytic number theory. For the series at $\rmi\infty$, the expansion reads
\begin{multline}
\label{eq:FourierExpEisInfty}
E_{k,p}^{\rmi\infty}(\tau;r_1,r_p)=1+(-1)^{\frac{k-k'}{2}}\frac{(2\uppi)^k}{\Gamma(k)p^k}\sum_{n \in \numgeq{Z}{1}}n^{k-1}\sum_{c \in \numgeq{Z}{1}}\frac{1}{c^k} \\
\times\sum_{\twoscript{0 \leq d < pc}{\gcd(pc, d)=1}}\etp{\frac{-n_\infty a+(n-n_\infty)d}{pc}-\frac{1}{2}\left(r_1s(-d,pc)+r_ps(-d,c)\right)}q^n,
\end{multline}
where $k'=(r_1+r_p)/2$, $n_\infty=(r_1+pr_p)/24$ and for each pair $(c,d)$, the letter $a$ stands for the inverse of $d$ mod $pc$. There is a similar expression for $E_{k,p}^{1}(\tau;r_1,r_p)$. See Theorem \ref{thm:FourierEisInfty} and Theorem \ref{thm:FourierEis1} for details.

As the third result, we explore when an Eisenstein series of rational weight may be equal to an eta-quotient. For level $2$ and $3$ we find out infinitely many such identities. For instance, for any rational number $n_1$ with $2/3 < n_1 \leq 1$ we have
\begin{equation*}
E_{3n_1,3}^{\rmi\infty}(\tau;9n_1,-3n_1)=\eta^{9n_1}(\tau)\eta^{-3n_1}(3\tau).
\end{equation*}
We also obtain identities for $E_{k,2}^{\rmi\infty}(\tau;r_1,r_2)$, $E_{k,3}^{1}(\tau;r_1,r_3)$ and $E_{k,2}^{1}(\tau;r_1,r_2)$. See Theorem \ref{thm:EisEquEtaInfty} and Theorem \ref{thm:EisEquEta1} for details.

Finally, we give an application which combines the above three results. We give series expressions of special values of Euler Gamma function at any rational arguments. For two typical examples, we have
\begin{align*}
\Gamma\left(\frac{8}{3}\right)&=-\frac{3}{32}\uppi^{8/3}\sum_{c\in\numgeq{Z}{1}}\frac{1}{c^{8/3}}\sum_{\twoscript{0 \leq d < 2c}{\gcd(2c,d)=1}}\etp{\frac{d}{2c}-\frac{8}{3}\left(2s(-d,2c)-s(-d,c)\right)}, \\
\Gamma\left(\frac{15}{7}\right)&=\frac{49}{540}\left(\frac{2\uppi}{3}\right)^{15/7}2^{8/7}\sum_{c\in\numgeq{Z}{1}}\frac{1}{c^{15/7}}\sum_{\twoscript{0 \leq d < 3c}{\gcd(3c,d)=1}}\etp{\frac{2d}{3c}-\frac{15}{14}\left(3s(-d,3c)-s(-d,c)\right)}.
\end{align*}
For the complete statement, see Theorem \ref{thm:GammaValue} and Corollary \ref{coro:GammaValue}. In the latter, real arguments are allowed.

Notations that are directly related to the above four results are collected in the beginning of Section \ref{sec:Holomorphic Eisenstein series of rational weights and prime levels}. All identities have been verified numerically by SageMath \cite{Sage} programs. The content of each section is implied by its section name, possibly except Section \ref{sec:Discussion}, where we discuss how and when Eisenstein series of rational weights and prime levels may be reduced to known functions.

\section{Modular forms of rational weights}
\label{sec:Modular forms of rational weights}
As one may use a double cover of the modular group when dealing with modular forms of half-integral weights (see \cite[\S 1]{Shi73}), it is necessary to use a multiple cover of the modular group for rational weights. Let $D$ be a positive integer. Then the $D$-cover of $\slR$ is defined as
\begin{equation*}
\Dcover{\slR}{D}=\left\{\elesltRaD{a}{b}{c}{d}{\varepsilon}{D} \colon \tbtmat{a}{b}{c}{d} \in \slR,\, \varepsilon^D=1\right\},
\end{equation*}
where the fractional power is the principal branch, that is, $z^\frac{m}{n}=\exp(\frac{m}{n}\log z)$ with $-\uppi < \Im(\log z) \leq \uppi$. The notation $\elesltRaDs{a}{b}{c}{d}{\varepsilon}{D}$ stands for $\elesltRaD{a}{b}{c}{d}{\varepsilon}{D}$ for simplicity\footnote{When $c=0$ and $d=-1$, one notation may be confused with the other which must be treated with caution.}. The group composition of $\Dcover{\slR}{D}$ is given by the following formula:
\begin{equation*}
\left(\tbtmat{a_1}{b_1}{c_1}{d_1}, \varepsilon_1\right)\left(\tbtmat{a_2}{b_2}{c_2}{d_2}, \varepsilon_2\right)=\left(\tbtmat{a_1}{b_1}{c_1}{d_1}\tbtmat{a_2}{b_2}{c_2}{d_2}, \varepsilon_1\varepsilon_2\delta\right),
\end{equation*}
where $\delta$ is a $D$-th root of unity determined by
\begin{equation*}
\delta=\frac{\sqrt[D]{c_1(a_2\tau+b_2)/(c_2\tau+d_2)+d_1}\sqrt[D]{c_2\tau+d_2}}{\sqrt[D]{c_1(a_2\tau+b_2)+d_1(c_2\tau+d_2)}}.
\end{equation*}
For a subgroup $G$, the notation $\Dcover{G}{D}$ means the preimage of $G$ under the natural projection $\Dcover{\slR}{D} \rightarrow \slR$.

Fix a positive integer $D$. Let $k \in \frac{1}{D}\numZ$, $G$ be a finite index subgroup of $\slZ$, and $\chi \colon \Dcover{G}{D}\rightarrow \numC^\times$ be a finite index character. The two prefixes ``finite index'' mean that $[\slZ \colon G] < +\infty$ and $[\Dcover{G}{D} \colon \ker \chi] < +\infty$. We define the slash operator of weight $k$ as
\begin{equation}
\label{eq:slash}
f\vert_k\elesltRaDs{a}{b}{c}{d}{\varepsilon}{D}(\tau)=\varepsilon^{-Dk}(c\tau+d)^{-k}f\left(\frac{a\tau+b}{c\tau+d}\right)
\end{equation}
for a continuous function $f\colon\uhp\rightarrow\numC$ and $\elesltRaDs{a}{b}{c}{d}{\varepsilon}{D} \in \Dcover{\slR}{D}$. It is a right group action\footnote{If we use the group $\slR$, instead of its $D$-cover, then \eqref{eq:slash} without the factor $\varepsilon^{-Dk}$ does not give a group action. This is one advantage of considering $D$-covers of modular groups. Another advantage is that the multiplier system is really a group character.} of $\Dcover{\slR}{D}$ on continuous functions (or holomorphic, or meromorphic functions). By a modular form of weight $k$, on the group $\Dcover{G}{D}$, and of character $\chi$, we mean a holomorphic function $f\colon \uhp \rightarrow \numC$ satisfying the following properties:
\begin{enumerate}
	\item $f\vert_k\gamma = \chi(\gamma)f$ for any $\gamma \in \Dcover{G}{D}$,
	\item $\lim_{\Im(\tau)\to+\infty}f\vert_k\gamma$ exists for any $\gamma \in \Dcover{\slZ}{D}$.
\end{enumerate}
The character $\chi$ is called the multiplier system of $f$. The complex vector space of all modular forms is denoted by $M_k(\Dcover{G}{D};\chi)$. The reader may compare our definition with Ibukiyama's (\cite{Ibu00} and \cite{Ibu20_2}). Ibukiyama constructs modular forms of rational weights using quotients of certain theta series and fractional power of the eta function. However, we think it is necessary and extremely convenient to use $D$-covers of modular groups when dealing with formal Eisenstein series of rational weights, as we shall do next.

\section{Eisenstein series of rational weights}
\label{sec:Eisenstein series of rational weights}
First some notations. If a group $G$ acts from the left (right resp.) on the set $X$, then we use $G\backslash X$ ($X/G$ resp.) to denote the set of orbits. If the subgroup $H$ acts on $G$, then the action is always understood as translations. For an element $g \in G$, $\langle g\rangle$ denotes the subgroup generated by $g$. The group considered in this section is $\Dcover{\Gamma_0(N)}{D}$ with $N, D$ positive integers. The symbol $\projQ$ denotes $\numQ$ together with the infinity point $\rmi\infty$ as a subset of the Riemann sphere. Notice $a/0=\rmi\infty$ for nonzero integer $a$. The group $\Gamma_0(N)$ and so $\Dcover{\Gamma_0(N)}{D}$ act (from the left) on $\projQ$ as $\tbtmat{a}{b}{c}{d}s=\frac{as+b}{cs+d}$. For $s \in \projQ$, $\overline{s}$ denotes its image in $\Dcover{\Gamma_0(N)}{D}\backslash\projQ$, and is called a cusp. The width of $\overline{s}$, denoted by $w_{\overline{s}}$, is defined as $[\pslZ_{s}\colon\overline{\Gamma_0(N)_s}]$, where the notation $G_s$ means the isotropy group, and $\overline{\Gamma_0(N)}$ is the quotient of $\Gamma_0(N)$ by $\{\pm I\}$. If $\gamma=\tbtmat{a}{b}{c}{d} \in \slR$, then the notation $\widetilde{\gamma}$ denotes $\elesltRaDs{a}{b}{c}{d}{1}{D}$ in $\Dcover{\slR}{D}$. Put $T=\tbtmat{1}{0}{0}{1}$, $S=\tbtmat{0}{-1}{1}{0}$ and $I=\tbtmat{1}{0}{0}{1}$.
\begin{deff}
\label{deff:Eis}
Let $D,N$ be positive integers, and let $k\in\frac{1}{D}\numZ$ that is greater than $2$. Let $\chi\colon \Dcover{\Gamma_0(N)}{D} \to \numC^\times$ be a finite index character such that $\chi(\widetilde{-I})=\etp{-k/2}$. Let $s \in \projQ$ and $\gamma_s \in \slZ$ such that $\gamma_s(\rmi\infty)=s$. If $\chi(\widetilde{\gamma_s}\widetilde{T}^{w_{\overline{s}}}\widetilde{\gamma_s}^{-1})=1$, then we define the Eisenstein series $E_{\gamma_s,k}$ on the group $\Dcover{\Gamma_0(N)}{D}$, of weight $k$, with character $\chi$ and at cusp $\gamma_s$ as
\begin{equation*}
E_{\gamma_s,k}(\tau)=\sum_{\gamma \in \widetilde{\gamma_s}\langle \widetilde{T}^{w_{\overline{s}}}\rangle\widetilde{\gamma_s}^{-1}\backslash\Dcover{\Gamma_0(N)}{D}}\chi(\gamma)^{-1}\cdot 1\vert_k\widetilde{\gamma_s}^{-1}\gamma.
\end{equation*}
\end{deff}
The condition $\chi(\widetilde{\gamma_s}\widetilde{T}^{w_{\overline{s}}}\widetilde{\gamma_s}^{-1})=1$ makes the formal series well-defined, that is, the terms are independent of the choice of representatives of $\widetilde{\gamma_s}\langle \widetilde{T}^{w_{\overline{s}}}\rangle\widetilde{\gamma_s}^{-1}\backslash\Dcover{\Gamma_0(N)}{D}$. The condition $k>2$ and $\chi$ is of finite index ensure that the series converges normally and hence defines a holomorphic function on $\uhp$.
\begin{lemm}
\label{lemm:EisAnoExpr}
Assumptions as in Definition \ref{deff:Eis}. Put $\gamma_s=\tbtmat{a_s}{b_s}{c_s}{d_s}$. Let $\gamma_1=\tbtmat{a_1}{b_1}{c_1}{d_1} \in \slZ$ be arbitrary. Then we have
\begin{equation*}
E_{\gamma_s,k}\vert_k\widetilde{\gamma_1}(\tau)=D\cdot\sum_{c,d}\chi\left(\widetilde{\tbtmat{a_s}{b_s}{c_s}{d_s}}\widetilde{\tbtmat{a}{b}{c}{d}}\widetilde{\tbtmat{a_1}{b_1}{c_1}{d_1}}^{-1}\right)^{-1}(c\tau+d)^{-k},
\end{equation*}
where the summation is over pairs of integers $(c,d)$ satisfying $\gcd(c,d)=1$ and $N \mid c_sad_1+d_scd_1-c_sbc_1-d_sdc_1$ for some $a,b$ with $ad-bc=1$.
\end{lemm}
\begin{proof}
By a change of variables, we may write
\begin{equation}
\label{eq:proofEisExpression}
E_{\gamma_s,k}\vert_k\widetilde{\gamma_1}=\sum_{\gamma \in \langle\widetilde{T}^{w_{\overline{s}}}\rangle\backslash\widetilde{\gamma_s}^{-1}\Dcover{\Gamma_0(N)}{D}\widetilde{\gamma_1}}\chi(\widetilde{\gamma_s}\gamma\widetilde{\gamma_1}^{-1})^{-1}\cdot 1\vert_k\gamma.
\end{equation}
It can be verified directly that the following set is a complete system of representatives of $\langle\widetilde{T}^{w_{\overline{s}}}\rangle\backslash\Dcover{\slZ}{D}$:
\begin{equation*}
\mathscr{R}=\left\{\widetilde{T}^r\cdot\elesltRaDs{a}{b}{c}{d}{\varepsilon}{D}\colon 0 \leq r \leq w_{\overline{s}}-1,\,c,d\in\numZ,\gcd(c,d)=1,\varepsilon^D=1\right\},
\end{equation*}
where for each pair $(c,d)$ in elements of $\mathscr{R}$, we shall fix a pair of integers $(a,b)$ with $ad-bc=1$. Therefore, the set
\begin{equation*}
\mathscr{R}_1=\{\gamma \in \mathscr{R}\colon \widetilde{\gamma_s}\gamma\widetilde{\gamma_1}^{-1}\in \Dcover{\Gamma_0(N)}{D}\}=\mathscr{R}\cap\widetilde{\gamma_s}^{-1}\Dcover{\Gamma_0(N)}{D}\widetilde{\gamma_1}
\end{equation*}
forms a complete system of representatives of $\langle\widetilde{T}^{w_{\overline{s}}}\rangle\backslash\widetilde{\gamma_s}^{-1}\Dcover{\Gamma_0(N)}{D}\widetilde{\gamma_1}$. Write out elements explicitly, we find that $\mathscr{R}_1$ consists of exactly those $\widetilde{T}^r\cdot\elesltRaDs{a}{b}{c}{d}{\varepsilon}{D}$ with $0 \leq r \leq w_{\overline{s}}-1$, $c,d\in\numZ,\gcd(c,d)=1$, $N \mid c_s(a+cr)d_1+d_scd_1-c_s(b+dr)c_1-d_sdc_1$ and $\varepsilon^D=1$. It can be verified that if both $\widetilde{T}^{r_1}\cdot\elesltRaDs{a}{b}{c}{d}{\varepsilon}{D}$ and $\widetilde{T}^{r_2}\cdot\elesltRaDs{a}{b}{c}{d}{\varepsilon}{D}$ are in $\mathscr{R}_1$, then $r_1=r_2$. Inserting this into \eqref{eq:proofEisExpression}, and taking into account the fact $\sum_{\varepsilon^D=1}\chi\left(I,\varepsilon\right)^{-1}\cdot \varepsilon^{-Dk}=D$ (since $\chi(\widetilde{-I})=\etp{-k/2}$ and $\widetilde{-I}^2=(I, \etp{1/D})$), the desired formula follows.
\end{proof}
\begin{thm}
\label{thm:EisMod}
Assumptions as in Definition \ref{deff:Eis}. Let $\gamma_1 \in \slZ$ be arbitrary. Then we have
\begin{equation*}
\lim_{\Im \tau \to +\infty}E_{\gamma_s,k}\vert_k\widetilde{\gamma_1}(\tau)=\begin{cases}
0,  & \text{if } \overline{\gamma_s(\rmi\infty)} \neq \overline{\gamma_1(\rmi\infty)}\\
2D\cdot\chi(\widetilde{g})^{-1}\cdot 1\vert_k\widetilde{\gamma_s}^{-1}\widetilde{g}\widetilde{\gamma_1}, & \text{else,} 
\end{cases}
\end{equation*}
where $g$ is any matrix in $\Gamma_0(N)$ satisfying $\gamma_s(\rmi\infty)=g\gamma_1(\rmi\infty)$. As a consequence,we have $E_{\gamma_s,k} \in M_k(\Dcover{\Gamma_0(N)}{D};\chi)$. Moreover, the set
\begin{equation*}
\left\{E_{\gamma_s,k} \colon \overline{s} \in \Dcover{\Gamma_0(N)}{D}\backslash\projQ,\, \chi(\widetilde{\gamma_s}\widetilde{T}^{w_{\widetilde{s}}}\widetilde{\gamma_s}^{-1})=1 \right\}
\end{equation*}
is $\numC$-linearly independent.
\end{thm}
\begin{proof}
Since the series \eqref{eq:proofEisExpression} converges uniformly on any compact subset of $\uhp$, we can change the order of the limit and the summation. Hence we obtain
\begin{equation*}
\lim_{\Im \tau \to +\infty}E_{\gamma_s,k}\vert_k\widetilde{\gamma_1}(\tau)=\sum_{\gamma}\chi(\widetilde{\gamma_s}\gamma\widetilde{\gamma_1}^{-1})^{-1}\cdot 1\vert_k\gamma,
\end{equation*}
where $\gamma$ ranges over $\langle\widetilde{T}^{w_{\overline{s}}}\rangle\backslash(\widetilde{\gamma_s}^{-1}\Dcover{\Gamma_0(N)}{D}\widetilde{\gamma_1}\cap\Dcover{\slZ}{D}_{\rmi\infty})$. If $\overline{\gamma_s(\rmi\infty)} \neq \overline{\gamma_1(\rmi\infty)}$, the summation range is empty, so the limit is zero. Otherwise, there is some $g \in \Gamma_0(N)$ such that $\widetilde{\gamma_s}^{-1}\widetilde{g}\widetilde{\gamma_1} \in \Dcover{\slZ}{D}_{\rmi\infty}$. Then the set $\widetilde{\gamma_s}^{-1}\Dcover{\Gamma_0(N)}{D}\widetilde{\gamma_1}\cap\Dcover{\slZ}{D}_{\rmi\infty}$ is a disjoint union of the the following $2D$ cosets:
\begin{equation*}
\langle\widetilde{T}^{w_{\overline{s}}}\rangle\cdot(\widetilde{-I})^t\widetilde{\gamma_s}^{-1}\widetilde{g}\widetilde{\gamma_1},\quad t=0,1,\dots,2D-1.
\end{equation*}
It follows that $\lim_{\Im \tau \to +\infty}E_{\gamma_s,k}\vert_k\widetilde{\gamma_1}(\tau)=2D\cdot\chi(\widetilde{g})^{-1}\cdot 1\vert_k\widetilde{\gamma_s}^{-1}\widetilde{g}\widetilde{\gamma_1}$, since $\sum_{t=0}^{2D-1}\chi(\widetilde{-I})^{-t}1\vert_k(\widetilde{-I})^t=2D$. The assertion on linear independency follows immediately. It remains to prove $E_{\gamma_s,k} \in M_k(\Dcover{\Gamma_0(N)}{D};\chi)$. The second condition of modular forms has been proved. To prove the first, let $\gamma_0 \in \Dcover{\Gamma_0(N)}{D}$ be arbitrary. We have, applying a change of variables,
\begin{equation*}
E_{\gamma_s,k}\vert_k\gamma_0=\chi(\gamma_0)\cdot\sum_{\gamma \in \widetilde{\gamma_s}\langle \widetilde{T}^{w_{\overline{s}}}\rangle\widetilde{\gamma_s}^{-1}\backslash\Dcover{\Gamma_0(N)}{D}}\chi(\gamma\gamma_0)^{-1}\cdot 1\vert_k\widetilde{\gamma_s}^{-1}\gamma\gamma_0=\chi(\gamma_0)E_{\gamma_s,k}.
\end{equation*}
This concludes the proof.
\end{proof}

\section{Eta quotients of rational weights}
\label{sec:Eta quotients of rational weights}
The simplest way to construct concrete examples of modular forms of rational weights, and to construct concrete characters of $D$-covers of modular groups is to use fractional powers of eta-quotients. For $r \in \numQ$, we define
\begin{align*}
\log\eta(\tau)&= \log\eta(\rmi)+\int_{\rmi}^{\tau}\frac{\eta'(z)}{\eta(z)}\diff z, \\
\eta^r(\tau)&=\exp\left(r\cdot\log\eta(\tau)\right).
\end{align*}
Note that $\eta(\rmi)$ is a real number, and $\log\eta(\rmi)$ is the real logarithm. The function $\log\eta$ is a single one, and should not be understood as the composition of $\log$ function and $\eta$ function. The propositions we will prove in the remaining section all rely on the classical formula of Dedekind:
\begin{equation*}
\log\eta\left(\frac{a\tau+b}{c\tau+d}\right)-\log\eta(\tau)=2\uppi\rmi\left(\frac{a+d}{24c}+\frac{1}{2}s(-d,c)-\frac{1}{8}\right)+\frac{1}{2}\log(c\tau+d),
\end{equation*}
where $\tbtmat{a}{b}{c}{d} \in \slZ$ with $c>0$, and $\log(c\tau+d)$ is the principal branch. For a proof, see \cite[Equation (12), \S 3.4]{Apo90}.

Let $N,D$ be positive integers. By an eta-quotient on $\Dcover{\Gamma_0(N)}{2D}$, we mean a product $\prod_{0<n \mid N}\eta(n\tau)^{r_n}$ with $r_n \in \frac{1}{D}\numZ$.
\begin{lemm}
\label{lemm:etaQuoChar}
Put $f(\tau)=\prod_{0<n \mid N}\eta(n\tau)^{r_n}$, and $k=\frac{1}{2}\sum_{n \mid N}r_n$. For any $\gamma=\elesltRaDs{a}{b}{c}{d}{\varepsilon}{2D} \in \Dcover{\Gamma_0(N)}{2D}$, we have $f\vert_k\gamma=\chi(\gamma)f$, where $\chi\colon\Dcover{\Gamma_0(N)}{2D} \to \numC^\times$ is a group character defined by $\chi(\gamma)=\varepsilon^{-2Dk}\etp{v(a,b,c,d)}$, where
\begin{equation*}
v(a,b,c,d)=\begin{cases}
\frac{a+d}{24c}\sum_{n \mid N}nr_n+\frac{1}{2}\sum_{n \mid N}r_ns(-d,\frac{c}{n})-\frac{k}{4} &\text{if } c>0 \\
\frac{a+d}{24c}\sum_{n \mid N}nr_n+\frac{1}{2}\sum_{n \mid N}r_ns(d,\frac{-c}{n})+\frac{k}{4} &\text{if } c<0 \\
\frac{b}{24}\sum_{n \mid N}nr_n &\text{if } c=0, a>0\\
-\frac{b}{24}\sum_{n \mid N}nr_n-\frac{k}{2} &\text{if } c=0, a<0.
\end{cases}
\end{equation*}
\end{lemm}
\begin{proof}
Note that $f(\tau)=\exp{\sum_{n \mid N}r_n\cdot\log\eta(n\tau)}$. Hence the assertion on the value of $\chi(\gamma)$ follows from the definition of slash operators, Dedekind's functional equation of $\log\eta$, the fact $\eta^{r_n}(n\tau)\vert_{r_n/2}\widetilde{T}=\etp{nr_n/24}\eta^{r_n}(n\tau)$, and the formula
\begin{equation*}
\elesltRaDs{-1}{0}{0}{-1}{1\cdot(-1)^{\frac{1}{2D}}}{2D}\cdot\elesltRaD{a}{b}{c}{d}{\varepsilon}{2D}= \elesltRaDs{-a}{-b}{-c}{-d}{\varepsilon\cdot(-c\tau-d)^{\frac{1}{2D}}}{2D} \text{    if } c<0.
\end{equation*}
The assertion that $\chi$ is a group character follows from the fact that slash operators give right group actions.
\end{proof}

If a function satisfies the first condition of the definition of modular forms, we say it transforms like a modular form. Hence an eta-quotient (of rational powers) transforms like a modular form. For a non-zero holomorphic function $f$ on $\uhp$ that transforms like a modular form and any $\gamma \in \Dcover{\slZ}{2D}$, $f\vert_k\gamma$ can be expressed as a function series\footnote{It can be proved that, as in the case of integral weight modular forms, this series converges normally for $\tau \in \uhp$.} of the form $\sum_{n \in \frac{1}{m}\numZ}a_{n}q^{n}$ where $m$ is a fixed positive integer. The least $n$ such that $a_n \neq 0$ is called the order of $f$ at cusp $\overline{\gamma(\rmi\infty)}$ (The order may be $-\infty$). We shall need the following formula concerning the order of an eta-quotient:
\begin{lemm}
\label{lemm:ordEtaQuo}
Notations and Assumptions as in Lemma \ref{lemm:etaQuoChar}. Let $a,c$ be coprime integers with $c>0$. Then the order of $f$ at the cusp $\overline{a/c} \in \Gamma_0(N)\backslash\projQ$ is $\frac{1}{24}\sum_{n \mid N}\frac{r_n\cdot\gcd(n,c)^2}{n}$.
\end{lemm}
\begin{proof}
If all $r_n$'s are integers, the desired formula is equivalent to \cite[Theorem 1.65]{Ono04} (Caution: the meaning of ``order'' at cusps here is different from that in \cite{Ono04}). Otherwise, apply the case just proven to $f^{D}$, and use the fact
\[\left(f\vert_k\elesltRaD{a}{b}{c}{d}{\varepsilon}{2D}\right)^{D_0}=f^{D_0}\vert_{D_0k}\elesltRaD{a}{b}{c}{d}{\varepsilon^{D_0}}{2D/D_0}\]
for $0 < D_0 \mid 2D$.
\end{proof}

\begin{rema}
\label{rema:etaCoeff}
The coefficient of the leading term of the $q$-expansion of any (rational-weight) eta-quotient must be $1$. This is not immediate since we have defined fractional powers of Dedekind eta function as the exponential of certain integral. We now prove this. Let $m$ be an integer and put $\eta^m(\tau)=q^{m/24}\sum_{n \in \numgeq{Z}{0}}a_nq^n$. Then obviously $a_0=1$.Put $\eta^{m/D}(\tau)=q^{m/(24D)}\sum_{n \in \numgeq{Z}{0}}b_nq^n$. Then $(\sum_{n \in \numgeq{Z}{0}}b_nq^n)^D=\sum_{n \in \numgeq{Z}{0}}a_nq^n$ and so $b_0^D=1$. We shall prove that $b_0=1$. Since the infinite product expansion of $\eta$ converges normally, we can take the logarithmic derivative and obtain a series expansion of $\eta'(z)/\eta(z)$ which still converges normally (see \cite[\S 1.2.3]{Rem98}). So we can change the order of the integral and the summation and obtain $\lim_{\Im\tau \to +\infty}(\log\eta(\tau)-\frac{2\uppi\rmi}{24}\tau) = 0$. Therefore
\begin{equation*}
b_0=\exp\lim_{\Im\tau \to +\infty}\frac{m}{D}\left(\log\eta(\tau)-\frac{2\uppi\rmi}{24}\tau\right)=1.
\end{equation*}
Consequently, for any eta-quotient $\prod_{n \mid N}\eta^{r_n}(n\tau)$ with $r_n$ rational numbers, the coefficient of $q^{\sum_{n \mid N}nr_n/24}$-term, which is the term with the least exponent, is $1$.
\end{rema}

We can construct modular forms of rational weights by eta-quotients of rational powers whose order at all cusps are non-negative. We present examples for level $11$ and multiplier systems as simple as possible in the remaining of this section. 

From a set of generators of $\Gamma_0(11)$ (can be obtained by a SageMath \cite{Sage} command \lstinline{Gamma0(11).generators()} which calculates generators using Farey symbols) we obtain a set of generators of $\Dcover{\Gamma_0(11)}{2D}$:
\begin{equation*}
\widetilde{T},\quad \widetilde{T}\widetilde{S}^{-1}\widetilde{T}^3\widetilde{S}^{-1}\widetilde{T}^4\widetilde{S},\quad \widetilde{T}\widetilde{S}^{-1}\widetilde{T}^4\widetilde{S}^{-1}\widetilde{T}^3\widetilde{S},\quad \widetilde{S}^2=\widetilde{-I}.
\end{equation*}
\begin{lemm}
\label{lemm:charGamma011}
Let $r_1$ and $r_{11}$ be rational numbers, and $D$ be a positive integer such that $D\cdot r_1$ and $D\cdot r_{11}$ are both integers. Then the multiplier system $\chi\colon\Dcover{\Gamma_0(11)}{2D}\to\numC^\times$ of the eta-quotient $\eta^{r_1}(\tau)\eta^{r_{11}}(11\tau)$ is determined by the following values:
\begin{align*}
\widetilde{T}&\mapsto \etp{\frac{r_1+11r_{11}}{24}} &\quad \widetilde{T}\widetilde{S}^{-1}\widetilde{T}^3\widetilde{S}^{-1}\widetilde{T}^4\widetilde{S}&\mapsto \etp{\frac{11r_1+13r_{11}}{24}}\\
\widetilde{S}^2&\mapsto\etp{\frac{-6r_1-6r_{11}}{24}} &\quad \widetilde{T}\widetilde{S}^{-1}\widetilde{T}^4\widetilde{S}^{-1}\widetilde{T}^3\widetilde{S}&\mapsto \etp{\frac{11r_1+13r_{11}}{24}}.
\end{align*}
\end{lemm}
\begin{proof}
Apply Lemma \ref{lemm:etaQuoChar} and take into account the fact
\begin{align*}
\widetilde{T}\widetilde{S}^{-1}\widetilde{T}^3\widetilde{S}^{-1}\widetilde{T}^4\widetilde{S}&=\elesltRaDs{7}{-2}{11}{-3}{\etp{-\frac{1}{2D}}(11\tau-3)^{\frac{1}{2D}}}{2D}, \\
\widetilde{T}\widetilde{S}^{-1}\widetilde{T}^4\widetilde{S}^{-1}\widetilde{T}^3\widetilde{S}&=\elesltRaDs{8}{-3}{11}{-4}{\etp{-\frac{1}{2D}}(11\tau-4)^{\frac{1}{2D}}}{2D}.
\end{align*}
\end{proof}

\begin{prop}
Among the functions $\eta^{r_1}(\tau)\eta^{r_{11}}(11\tau)$ with $r_1$ and $r_{11}$ rationals, the following ones are exactly those modular forms on $\Dcover{\Gamma_0(11)}{2D}$ for some positive integer $D$ such that the multiplier system is trivial on the subgroup generated by $\widetilde{T}$, $\widetilde{T}\widetilde{S}^{-1}\widetilde{T}^3\widetilde{S}^{-1}\widetilde{T}^4\widetilde{S}$ and $\widetilde{T}\widetilde{S}^{-1}\widetilde{T}^4\widetilde{S}^{-1}\widetilde{T}^3\widetilde{S}$:
\begin{equation*}
\eta(\tau)^{-26m_1-12m_2}\cdot\eta(11\tau)^{22m_1+12m_2}, \quad m_1\in\frac{1}{9}\numZ,\,m_2\in\numZ,\,-\frac{9}{5}m_1 \leq m_2 \leq -\frac{11}{5}m_1.
\end{equation*}
Moreover, one can choose $2D=18$.
\end{prop}
\begin{proof}
A consequence of Lemma \ref{lemm:ordEtaQuo} and Lemma \ref{lemm:charGamma011}.
\end{proof}

\section{A valence formula for rational weights}
\label{sec:A valence formula for rational weights}
We shall use the following fundamental fact on meromorphic modular forms of rational weights.
\begin{thm}
\label{thm:valence}
Let $D$ be a positive integer and $k \in \frac{1}{D}\numZ$. Let $G$ be a subgroup of $\slZ$ of finite index and $\Dcover{G}{D}$ denote its $D$-cover. Let $\chi\colon \Dcover{G}{D}\to\numC^\times$ be a linear character of finite index. Then for any nonzero meromorphic modular form $f$ on the group $\Dcover{G}{D}$, of weight $k$ and with multiplier system $\chi$, we have
\begin{equation*}
\sum_{\overline{\gamma} \in \overline{G}\backslash\pslZ}\ord\nolimits_{\rmi\infty}(f\vert_k \gamma) + \sum_{\tau \in \overline{G}\backslash \uhp}\frac{1}{\abs{\overline{G}_\tau}}\ord\nolimits_\tau(f)=\frac{1}{12}[\pslZ : \overline{G}]k.
\end{equation*}
\end{thm}
By a meromorphic modular form, we mean a function that satisfies modular transformation laws (that is, transforms like a modular form), and is meromorphic on $\uhp$ and at all cusps of $G\backslash \projQ$ (which means the order of this function at any cusp is finite).

This formula is well-known, classical for $D=1,2$ and is called the valence formula. A detailed proof for $D=2$ is given in \cite[Theorem 2.1]{ZZ21}. The proof for general $D$ is almost the same as that for $D=2$ and we shall describe several slight changes now. Firstly, all double covers of modular groups in the proof of \cite[Theorem 2.1]{ZZ21} should be replaced by $D$-covers. Secondly, since the group of linear characters of the double cover $\Dcover{\slZ}{2}$ is a cyclic group of order $24$, we consider the $24$-th power of $\prod_if\vert_k\widetilde{M_i}$ to obtain a meromorphic modular form on the full modular group with trivial multiplier system. However, the group of linear characters of $\Dcover{\slZ}{D}$ has $12D$ elements (which we will prove in the next lemma). Therefore, we shall apply the usual valence formula to the $12D$-th power of $\prod_if\vert_k\widetilde{M_i}$, instead of the $24$-th power. The proof of \cite[Theorem 2.1]{ZZ21} with the above changes, together with the following lemma now serves a proof of Theorem \ref{thm:valence}.

\begin{lemm}
Put $T=\tbtmat{1}{1}{0}{1}$, $S=\tbtmat{0}{-1}{1}{0}$, and $I=\tbtmat{1}{0}{0}{1}$. Then the group $\Dcover{\slZ}{D}$ has a presentation with generators $\widetilde{T}$, $\widetilde{S}$ and relations $\widetilde{S}^{4D}=\widetilde{I}$, $(\widetilde{S}\widetilde{T})^3=\widetilde{S}^2$, $\widetilde{S}^4\widetilde{T}=\widetilde{T}\widetilde{S}^4$. As a consequence, there are exactly $12D$ linear characters on $\Dcover{\slZ}{D}$ which form a cyclic group.
\end{lemm}
\begin{proof}
Since $\widetilde{S}^4=\elesltRaDs{1}{0}{0}{1}{\etp{1/D}\cdot(0\tau+1)^{1/D}}{D}$ and $\slZ$ is generated by $T$ and $S$, we deduce that $\Dcover{\slZ}{D}$ is generated by $\widetilde{T}$ and $\widetilde{S}$. The three relations can be verified by a direct calculation. To prove that these generate all relations, we use more formal language to make the assertion precise. Let $x_T, x_S$ be two formal symbols. Put $X=\{x_T,x_S\}$. Let $F(X)$ be the free group generated by $X$. Let $p\colon F(X)\to \Dcover{\slZ}{D}$ be the group homomorphism that sends $x_T$ to $\widetilde{T}$ and $x_S$ to $\widetilde{S}$. Similarly, let $p_0\colon F(X)\to \slZ$ be the homomorphism that sends $x_T$, $x_S$ to $T$, $S$ respectively. Presentations of $\slZ$ are well-known; one of these gives that $\ker(p_0)=\langle x_S^4, x_S^2(x_Sx_T)^{-3}\rangle_{n}$ and $p_0$ is surjective, where $\langle H\rangle_n$ means the \emph{normal} subgroup generated by the subset $H$. Put
\begin{equation*}
K=\langle x_S^2(x_Sx_T)^{-3}, x_S^{4D}, x_T^{-1}x_S^{-4}x_Tx_S^4\rangle_n.
\end{equation*}
It remains to prove that $p$ induces an isomorphism from $F(X)/K$ onto $\Dcover{\slZ}{D}$. We have already proved that $p$ is surjective and $K \subseteq \ker(p)$. So let $r \in \ker(p)$ be arbitrary and we need to prove $r\cdot K=K$. Since $\ker(p)\subseteq\ker(p_0)$, we can write
\begin{equation*}
r=\prod_{i=1}^{t}g_iy_ig_i^{-1},
\end{equation*}
where each $y_i$ is $x_S^4$, $x_S^2(x_Sx_T)^{-3}$ or their inverses, and $g_i \in F(X)$. Since $x_T^{-1}x_S^{-4}x_Tx_S^4 \in K$, there is some $e \in \numZ$ such that
\begin{equation*}
\prod_{i=1}^{t}g_iy_ig_i^{-1}\cdot K=\left(\prod_{i=1}^{t'}h_iz_ih_i^{-1}\right)x_S^{4e}\cdot K,
\end{equation*}
where each $z_i$ is $x_S^2(x_Sx_T)^{-3}$ or its inverse, and $h_i \in F(X)$. Put $r'=\left(\prod_{i=1}^{t'}h_iz_ih_i^{-1}\right)x_S^{4e}$. It follows from the fact $K \subseteq \ker(p)$ and $r \in \ker(p)$ that $r' \in \ker(p)$. Hence $\left(\prod_{i=1}^{t'}p(h_i)p(z_i)p(h_i)^{-1}\right)\widetilde{S}^{4e}=\widetilde{I}$, and so $\widetilde{S}^{4e}=\widetilde{I}$. Therefore $D \mid e$, and $r\cdot K=r'\cdot K=K$. Finally, from the presentation, the linear characters on $\Dcover{\slZ}{D}$ are exactly the followings:
\begin{equation*}
\widetilde{T}\mapsto \etp{\frac{i}{12D}},\qquad \widetilde{S}\mapsto \etp{-\frac{i}{4D}},
\end{equation*}
where $i$ ranges over $1,2,\dots,12D$. These characters obviously form a cyclic group.
\end{proof}

What we need in the subsequent discussion is exactly the following corollary.
\begin{coro}
\label{coro:noCuspForm}
Let $p$ be a prime and $r_1, r_p$ be two rationals. Put $n_\infty=(r_1+pr_p)/24$, $n_1=(pr_1+r_p)/24$ and $k'=(r_1+r_p)/2$. Let $D,m_\infty,m_1$ be positive integers such that $Dr_1$, $Dr_p$, $m_\infty n_\infty$ and $m_1n_1$ are integers. Let $\chi_{r_1,r_p}\colon \Dcover{\Gamma_0(p)}{2D}\to \numC^\times$ be the multiplier system of $\eta^{r_1}(\tau)\eta^{r_p}(p\tau)$. If the following inequality is fulfilled:
\begin{equation*}
\frac{1}{m_\infty}+\frac{1}{m_1} > n_\infty+n_1,
\end{equation*}
then the space $M_{k'}(\Dcover{\Gamma_0(p)}{2D};\chi_{r_1,r_p})$ contains no nonzero cusp form.
\end{coro}
By a cusp form, we mean a modular form whose order at any cusp is positive.
\begin{proof}
We give an indirect proof. Suppose $f$ is a nonzero cusp form in $M_{k'}(\Dcover{\Gamma_0(p)}{2D};\chi_{r_1,r_p})$. Since $[\pslZ \colon \overline{\Gamma_0(p)}]=p+1$, we deduce from Theorem \ref{thm:valence} that
\begin{equation*}
\sum_{\overline{\gamma} \in \overline{G}\backslash\pslZ}\ord\nolimits_{\rmi\infty}(f\vert_{k'} \gamma) \leq n_\infty+n_1.
\end{equation*}
The quotient $\Gamma_0(p)\backslash\projQ$ contains exactly two cusps: $\overline{1/1}$ and $\overline{1/p}=\overline{\rmi\infty}$, whose widths are $p$ and $1$ respectively. Put $\gamma_1=\tbtmat{1}{0}{1}{1}$. Then we have $\ord\nolimits_{\rmi\infty}(f) + p\cdot\ord\nolimits_{\rmi\infty}(f\vert_{k'}\widetilde{\gamma_1}) \leq n_\infty+n_1$. We shall prove that $\ord\nolimits_{\rmi\infty}(f) \geq 1/m_\infty$ and $\ord\nolimits_{\rmi\infty}(f\vert_{k'}\widetilde{\gamma_1}) \geq 1/(pm_1)$, and hence reach a contradiction. To prove these two inequalities, it is necessary to calculate $\chi_{r_1,r_p}(\widetilde{T})$ and $\chi_{r_1,r_p}(\widetilde{\gamma_1}\widetilde{T}^p\widetilde{\gamma_1}^{-1})$ first. It is obviously that $\chi_{r_1,r_p}(\widetilde{T})=\etp{n_\infty}$. For $\chi_{r_1,r_p}(\widetilde{\gamma_1}\widetilde{T}^p\widetilde{\gamma_1}^{-1})$, using Lemma \ref{lemm:etaQuoChar} and the identity
\begin{equation*}
\widetilde{\tbtMat{1}{0}{1}{1}}\widetilde{\tbtMat{1}{p}{0}{1}}\widetilde{\tbtMat{1}{0}{1}{1}}^{-1}=\widetilde{\tbtMat{1-p}{p}{-p}{p+1}},
\end{equation*}
we obtain that
\begin{equation*}
\chi_{r_1,r_p}(\widetilde{\gamma_1}\widetilde{T}^p\widetilde{\gamma_1}^{-1})=\etp{\frac{r_1}{24}\left(-\frac{2}{p}+12s(1,p)+3\right)+\frac{r_p}{24}}.
\end{equation*}
By the reciprocity law for Dedekind sums (see \cite[Theorem 3]{RW41} for instance), we have $s(1,p)=(p-1)(p-2)/(12p)$. Hence $\chi_{r_1,r_p}(\widetilde{\gamma_1}\widetilde{T}^p\widetilde{\gamma_1}^{-1})=\etp{n_1}$. Now we have $f\vert_{k'}\widetilde{T}=\etp{n_\infty}f$, so $f$ has an expansion $\sum_{n \in \numZ}a_nq^{n/m_\infty}$. Since $f$ is a cusp form, we obtain $\ord\nolimits_{\rmi\infty}(f) \geq 1/m_\infty$. For another cusp, note that $f\vert_{k'}\widetilde{\gamma_1}\vert_{k'}\widetilde{T}^p=\etp{n_1}f\vert_{k'}\widetilde{\gamma_1}$, which implies that $f\vert_{k'}\widetilde{\gamma_1}$ has an expansion of the form $\sum_{n \in \numZ}b_nq^{n/(pm_1)}$. Therefore $\ord\nolimits_{\rmi\infty}(f\vert_{k'}\widetilde{\gamma_1}) \geq 1/(pm_1)$. This concludes the indirect proof.
\end{proof}

\section{Eisenstein series of rational weights and prime levels}
\label{sec:Holomorphic Eisenstein series of rational weights and prime levels}
Throughout this and remaining sections, we adopt the notations in Corollary \ref{coro:noCuspForm}: $p$ is a prime, $r_1, r_p$ are rational numbers, $n_\infty=(r_1+pr_p)/24$, $n_1=(pr_1+r_p)/24$, $k'=(r_1+r_p)/2$, $D, m_\infty, m_1$ are positive integers such that $Dr_1$, $Dr_p$, $m_\infty n_\infty$ and $m_1n_1$ are integers. Let $k \in{k'+2\numZ}$ be a rational number greater than $2$. Moreover, the character $\chi_{r_1,r_p}\colon \Dcover{\Gamma_0(p)}{2D}\to \numC^\times$ is the multiplier system of $\eta^{r_1}(\tau)\eta^{r_p}(p\tau)$ and $T=\tbtmat{1}{1}{0}{1}$, $\gamma_1=\tbtmat{1}{0}{1}{1}$. The following is our (first) main result.
\begin{thm}
\label{thm:mainEis}
The condition $\chi_{r_1,r_p}(\widetilde{T})=1$ holds if and only if $n_\infty \in \numZ$. In this case the series $E_{k,p}^{\rmi\infty}(\tau;r_1,r_p)$ (defined in \eqref{eq:Eisinfty}) belongs to the space $M_{k}(\Dcover{\Gamma_0(p)}{2D};\chi_{r_1,r_p})$. On the other hand, the condition $\chi_{r_1,r_p}(\widetilde{\gamma_1}\widetilde{T}^p\widetilde{\gamma_1}^{-1})=1$ holds if and only if $n_1 \in \numZ$. In this case we can define
\begin{multline*}
\label{eq:Eis1}
E_{k,p}^{1}(\tau;r_1,r_p)=\sum_{\twoscript{c > 0,d\in\numZ}{\gcd(c,d)=\gcd(p,c)=1}}(c\tau+d)^{-k} \times \\
\etp{-r_1\left(\frac{a+b+d}{24(a+c)}+\frac{1}{2}s(-b-d,a+c)-\frac{1}{8}\right)-r_p\left(\frac{a+b+d}{24(a+c)/p}+\frac{1}{2}s(-b-d,(a+c)/p)-\frac{1}{8}\right)},
\end{multline*}
where for each pair $(c,d)$, $a$ and $b$ are any two integers satisfying $ad-bc=1$, $p \mid a+c$ and $a+c>0$. This series belongs to $M_{k}(\Dcover{\Gamma_0(p)}{2D};\chi_{r_1,r_p})$.
\end{thm}
\begin{proof}
We have proved that $\chi_{r_1,r_p}(\widetilde{T})=\etp{n_\infty}$ and $\chi_{r_1,r_p}(\widetilde{\gamma_1}\widetilde{T}^p\widetilde{\gamma_1}^{-1})=\etp{n_1}$ in the proof of Corollary \ref{coro:noCuspForm}. The two ``if and only if'' statements follow from this. Let $E_{I,k}$ and $E_{\gamma_1,k}$ denote Eisenstein series in Definition \ref{deff:Eis}. (The level $N$ in Definition \ref{deff:Eis} should be replaced by $p$ here and $D$ in Definition \ref{deff:Eis} replaced by $2D$.) It has been proved that if $\chi_{r_1,r_p}(\widetilde{T})=1$ then $E_{I,k} \in M_{k}(\Dcover{\Gamma_0(p)}{2D};\chi_{r_1,r_p})$, and if $\chi_{r_1,r_p}(\widetilde{\gamma_1}\widetilde{T}^p\widetilde{\gamma_1}^{-1})=1$ then $E_{\gamma_1,k} \in M_{k}(\Dcover{\Gamma_0(p)}{2D};\chi_{r_1,r_p})$ in Theorem \ref{thm:EisMod}. We may conclude our proof if we can show $E_{I,k}(\tau)=4D\cdot E_{k,p}^{\rmi\infty}(\tau;r_1,r_p)$ and $E_{\gamma_1,k}(\tau)=4D\cdot E_{k,p}^{1}(\tau;r_1,r_p)$ under respective conditions. Since the proof for these two identities are similar, we only prove the latter (under the condition $n_1 \in \numZ$). Put $\chi=\chi_{r_1, r_p}$ for simplicity. By Lemma \ref{lemm:EisAnoExpr} we have
\begin{equation*}
E_{\gamma_1,k}(\tau)=2D\cdot\sum_{c,d}\chi\left(\widetilde{\tbtmat{1}{0}{1}{1}}\widetilde{\tbtmat{a}{b}{c}{d}}\right)^{-1}(c\tau+d)^{-k},
\end{equation*}
where the summation is over coprime integers $c, d$ such that there exist integers $a$ and $b$ with $ad-bc=1$ and $p \mid a+c$. Using elementary number theory we know that such pairs $(c,d)$ are exactly those with the property $\gcd(c,p)=\gcd(c,d)=1$. Moreover, a direct calculation shows that if $c<0$ then
\begin{equation*}
\chi\left(\widetilde{\tbtmat{1}{0}{1}{1}}\widetilde{\tbtmat{a}{b}{c}{d}}\right)^{-1}(c\tau+d)^{-k}=\chi\left(\widetilde{\tbtmat{1}{0}{1}{1}}\widetilde{\tbtmat{-a}{-b}{-c}{-d}}\right)^{-1}(-c\tau-d)^{-k}.
\end{equation*}
Here we have used the fact $\chi(\widetilde{-I})=\etp{-k/2}$ and
\begin{equation*}
((-1)(-c\tau-d))^{-k}=\etp{k}(-1)^{-k}(-c\tau-d)^{-k} \text{ if } c<0.
\end{equation*}
Therefore
\begin{equation}
\label{eq:E1kFormula}
E_{\gamma_1,k}(\tau)=4D\cdot\sum_{\twoscript{c > 0,d\in\numZ}{\gcd(c,d)=\gcd(p,c)=1}}\chi\left(\widetilde{\tbtmat{1}{0}{1}{1}}\widetilde{\tbtmat{a}{b}{c}{d}}\right)^{-1}(c\tau+d)^{-k}.
\end{equation}
Since each term in the above series is independent of the choice of $a$ and $b$, we can assume $a+c>0$. Hence
\begin{equation*}
\widetilde{\tbtMat{1}{0}{1}{1}}\widetilde{\tbtMat{a}{b}{c}{d}}=\widetilde{\tbtMat{a}{b}{a+c}{b+d}}.
\end{equation*}
Inserting Lemma \ref{lemm:etaQuoChar} into \eqref{eq:E1kFormula} gives the desired identity.
\end{proof}
We now turn to the aim of computing Fourier coefficients. Since there are some differences between methods for dealing with $E_{k,p}^{\rmi\infty}(\tau;r_1,r_p)$ and $E_{k,p}^{1}(\tau;r_1,r_p)$, we state them separately. First we consider $E_{k,p}^{\rmi\infty}(\tau;r_1,r_p)$.
\begin{thm}
\label{thm:FourierEisInfty}
Notations as in Theorem \ref{thm:mainEis}. Assume that $n_\infty$ is an integer. Then the Fourier expansion of $E_{k,p}^{\rmi\infty}(\tau;r_1,r_p)$ is given by the expression \eqref{eq:FourierExpEisInfty}.
\end{thm}
\begin{proof}
We introduce an auxiliary function $f(d,c)$ as follows:
\begin{equation*}
f(d,c)= \etp{-r_1\left(\frac{a+d}{24c}+\frac{1}{2}s(-d,c)\right)-r_p\left(\frac{a+d}{24c/p}+\frac{1}{2}s(-d,c/p)\right)},
\end{equation*}
where $p\mid c > 0$ and $d \in \numZ$ with $\gcd(c,d)=1$. Thus $E_{k,p}^{\rmi\infty}(\tau;r_1,r_p)=1+\etp{k'/4}\sum_{c,d}f(d,c)(c\tau+d)^{-k}$. We extend the domain of definition of $f$: if $\gcd(c,d) = 1$, $c \in \numgeq{Z}{1}$, $d \in \numZ$ but $p \nmid c$, then $f(d,c)=0$, and for the general case
\begin{equation*}
f(d,c)=f\left(\frac{d}{\gcd(c,d)},\frac{c}{\gcd(c,d)}\right).
\end{equation*}
Note that $f(d,c)=f(d+nc,c)$ for $n \in \numZ$ since $n_\infty \in \numZ$ and $s(-d,c)=s(-d+nc,c)$. From the fact $\sum_{t \mid \gcd(c,d)}\mu(t) = 0$ if $\gcd(c,d) > 1$, we obtain
\begin{equation}
\label{eq:proofFourierEisInfty}
\sum_{\twoscript{p \mid c > 0, d \in \numZ}{\gcd(c,d)=1}}f(d,c)(c\tau+d)^{-k}=\sum_{t \in \numgeq{Z}{1}}\frac{\mu(t)}{t^k}\sum_{c \in \numgeq{Z}{1}, d \in \numZ}f(d,c)(c\tau+d)^{-k},
\end{equation}
where $\mu$ denotes the M\"obius function. Using the Lipschitz summation formula
\begin{equation*}
\sum_{n \in \numZ}\frac{1}{(\tau+n)^s}=\rme^{-\uppi\rmi s/2}\frac{(2\uppi)^s}{\Gamma(s)}\sum_{n \geq 1}n^{s-1}\rme^{2\uppi\rmi n\tau},\quad \tau\in\uhp,\,\Re(s)>1,
\end{equation*}
we deduce that
\begin{align*}
\sum_{d \in \numZ}f(d,c)(c\tau+d)^{-k} &= c^{-k}\sum_{0 \leq d_1 < c}f(d_1,c)\sum_{d_0 \in \numZ}(\tau+\frac{d_1}{c}+d_0)^{-k} \\
&=\rme^{-\uppi\rmi k/2}\frac{(2\uppi)^k}{\Gamma(k)}c^{-k}\sum_{n \in \numgeq{Z}{1}}n^{k-1}\sum_{0 \leq d < c}f(d,c)\etp{\frac{dn}{c}}q^n.
\end{align*}
Inserting this into $\eqref{eq:proofFourierEisInfty}$, and using $\sum_{t}\frac{\mu(t)}{t^k}=\zeta(k)^{-1}$, we obtain
\begin{equation}
\label{eq:eq:proofFourierEisInfty2}
E_{k,p}^{\rmi\infty}(\tau;r_1,r_p)=1+(-1)^{\frac{k-k'}{2}}\frac{(2\uppi)^k}{\Gamma(k)\zeta(k)}\sum_{n \in \numgeq{Z}{1}}n^{k-1}\sum_{c \in \numgeq{Z}{1}}\frac{1}{c^k}\sum_{0 \leq d < c}f(d,c)\etp{\frac{dn}{c}}q^n.
\end{equation}
Regrouping the terms according to $\gcd(c,d)$, we have
\begin{multline*}
\sum_{c \in \numgeq{Z}{1}}\frac{1}{c^k}\sum_{0 \leq d < c}f(d,c)\etp{\frac{dn}{c}}=\frac{\zeta(k)}{p^k}\sum_{c\in\numgeq{Z}{1}}\frac{1}{c^k}\times \\
\sum_{\twoscript{0 \leq d < pc}{\gcd(pc,d)=1}}\etp{\frac{-n_\infty a+(n-n_\infty)d}{pc}-\frac{1}{2}\left(r_1s(-d,pc)+r_ps(-d,c)\right)},
\end{multline*}
where $a$ is an integer with $ad \equiv 1 \bmod pc$ for each pair $(c,d)$. The desired identity follows from this and \eqref{eq:eq:proofFourierEisInfty2}.
\end{proof}

Now we consider $E_{k,p}^{1}(\tau;r_1,r_p)$.
\begin{thm}
\label{thm:FourierEis1}
Notations as in Theorem \ref{thm:mainEis}. Assume that $n_1$ is an integer. Then the Fourier expansion of $E_{k,p}^{1}(\tau;r_1,r_p)$ is given by the following formula
\begin{multline*}
E_{k,p}^{1}(\tau;r_1,r_p)=(-1)^{\frac{k-k'}{2}}\frac{(2\uppi)^k}{\Gamma(k)m_\infty^k}\sum_{n \in \numgeq{Z}{1}}n^{k-1}\sum_{\twoscript{c \in \numgeq{Z}{1}}{\gcd(p,c)=1}}\frac{1}{c^k}\sum_{\twoscript{0 \leq d < cm_\infty}{\gcd(c, d)=1}} \\
\etp{\frac{dn}{cm_\infty}-n_\infty\frac{a+b+d}{a+c}-\frac{1}{2}\left(r_1s(-b-d,a+c)+r_ps(-b-d,(a+c)/p)\right)}q^{n/m_\infty},
\end{multline*}
where for each pair $(c,d)$, $a$ and $b$ are any two integers satisfying $ad-bc=1$, $p \mid a+c$ and $a+c>0$.
\end{thm}
\begin{proof}
As in the proof of Theorem \ref{thm:FourierEisInfty}, we introduce an auxiliary function $g(d,c)$ as follows:
\begin{equation*}
g(d,c)= \etp{-r_1\left(\frac{a+b+d}{24(a+c)}+\frac{1}{2}s(-b-d,a+c)\right)-r_p\left(\frac{a+b+d}{24(a+c)/p}+\frac{1}{2}s(-b-d,(a+c)/p)\right)},
\end{equation*}
where $a$ and $b$ are any two integers satisfying $ad-bc=1$, $p \mid a+c$ and $a+c>0$. It is defined for $c \in \numgeq{Z}{1}$, $d \in \numZ$ satisfying $\gcd(p,c)=\gcd(c,d)=1$. Thus $E_{k,p}^{1}(\tau;r_1,r_p)=\etp{k'/4}\sum_{c,d}g(d,c)(c\tau+d)^{-k}$. Extend the domain of definition as follows. If $c \in \numgeq{Z}{1}$, $d \in \numZ$, $\gcd(c,d)=1$ but $p \mid c$, then $g(d,c)=0$. If in general $c \in \numgeq{Z}{1}$ and $d \in \numZ$, then set
\begin{equation*}
g(d,c)=g\left(\frac{d}{\gcd(c,d)},\frac{c}{\gcd(c,d)}\right).
\end{equation*}
Using a similar argument as in the deduction of \eqref{eq:proofFourierEisInfty}, we obtain
\begin{equation*}
E_{k,p}^{1}(\tau;r_1,r_p)=\zeta(k)^{-1}\etp{k'/4}\sum_{c \in \numgeq{Z}{1}, d \in \numZ}g(d,c)(c\tau+d)^{-k}.
\end{equation*}
If we write out the occurrence of $a$ and $b$, we obtain
\begin{equation*}
g(d,c)=g(d,c;a,b)=g(d+cm_\infty,c;a,b+am_\infty)=g(d+cm_\infty,c),
\end{equation*}
which means for fixed $c$ the function $g(d,c)$ of $d$ is periodic with a period $cm_\infty$. Now an argument similar to the one used to deduce \eqref{eq:eq:proofFourierEisInfty2} shows that
\begin{equation}
\label{eq:proofFourierEis1}
E_{k,p}^{1}(\tau;r_1,r_p)=(-1)^{\frac{k-k'}{2}}\frac{(2\uppi)^k}{\Gamma(k)\zeta(k)m_\infty^{k}}\sum_{n \in \numgeq{Z}{1}}n^{k-1}\sum_{c \in \numgeq{Z}{1}}\frac{1}{c^k}\sum_{0 \leq d < cm_\infty}g(d,c)\etp{\frac{dn}{cm_\infty}}q^{n/m_\infty}.
\end{equation}
Regrouping the terms according to $\gcd(c,d)$, and using the fact $g(td,tc)=g(d,c)$, we obtain
\begin{multline*}
\sum_{c \in \numgeq{Z}{1}}\frac{1}{c^k}\sum_{0 \leq d < cm_\infty}g(d,c)\etp{\frac{dn}{cm_\infty}}=\zeta(k)\sum_{\twoscript{c \in \numgeq{Z}{1}}{\gcd(p,c)=1}}\frac{1}{c^k}\sum_{\twoscript{0 \leq d < cm_\infty}{\gcd(c, d)=1}} \\
\etp{\frac{dn}{cm_\infty}-n_\infty\frac{a+b+d}{a+c}-\frac{1}{2}\left(r_1s(-b-d,a+c)+r_ps(-b-d,(a+c)/p)\right)},
\end{multline*}
where $a$ and $b$ are as before. The desired identity follows from this and \eqref{eq:proofFourierEis1}.
\end{proof}

\begin{examp}
Set $p=11$, $r_1=44/9$, $r_{11}=-4/9$ and $k=20/9$. Then we have
\begin{equation*}
E_{20/9,11}^{\rmi\infty}(\tau;44/9,-4/9)=1+\etp{\frac{5}{9}}\sum_{\twoscript{11 \mid c > 0,\,d \in \numZ}{\gcd(c,d)=1}}\etp{-\frac{22}{9}s(-d, c)+\frac{2}{9}s(-d,c/11)}(c\tau+d)^{-20/9},
\end{equation*}
which is an Eisenstein series in the space $M_{20/9}(\Dcover{\Gamma_0(11)}{18},\chi)$. The character $\chi\colon\Dcover{\Gamma_0(11)}{18}\to\numC^\times$ is the one in Lemma \ref{lemm:charGamma011} with $r_1=44/9$, $r_{11}=-4/9$. The Fourier expansion reads
\begin{multline*}
E_{20/9,11}^{\rmi\infty}(\tau;44/9,-4/9)=1+\frac{(2\uppi)^{20/9}}{\Gamma(20/9)11^{20/9}}\\
\times\sum_{n \in \numgeq{Z}{1}}n^{11/9}\sum_{c \in \numgeq{Z}{1}}\frac{1}{c^{20/9}} \sum_{\twoscript{0 \leq d < 11c}{\gcd(11c, d)=1}}\etp{\frac{nd}{11c}}\etp{-\frac{2}{9}\left(11s(-d,11c)-s(-d,c)\right)}q^n.
\end{multline*}
\end{examp}

We will show in Section \ref{sec:Discussion} how $E_{k,p}^{\rmi\infty}(\tau;r_1,r_p)$ and $E_{k,p}^{1}(\tau;r_1,r_p)$ may be reduced to known functions when $r_1$ and $r_p$ are integers. In general cases, we can simplify neither exponential sums involving Dedekind sums nor $L$-values of these exponential sums at $k$ in Fourier expansions of $E_{k,p}^{\rmi\infty}(\tau;r_1,r_p)$ and $E_{k,p}^{1}(\tau;r_1,r_p)$. But in some special cases, we can evaluate the $L$-values at $k$ for any fixed $n$ using dimension arguments, which will be done in the next section.

\section{When an Eisenstein series equals an eta-quotient}
\label{sec:When an Eisenstein series equals an eta-quotient}
If a space of modular forms (with weight greater than $2$) is one-dimensional, and it contains both an eta-quotient (of rational exponents) and an Eisenstein series coincidentally, then they must be proportional. More generally, if a space of modular forms contains no nonzero cusp form, and contains an eta-quotient that can cancel constant terms of an Eisenstein series at both cusps, then they are equal. We find out all\footnote{If one uses stronger estimation of dimensions of spaces of cusp forms, then there may be more such identities.} possible identities of this type where dimensions are estimated by Theorem \ref{thm:valence} or its corollaries.

For $E_{k,p}^{\rmi\infty}(\tau;r_1,r_p)$, we have $m_\infty=1$. Moreover, if $E_{k,p}^{\rmi\infty}(\tau;r_1,r_p)=\eta^{r_1}(\tau)\eta^{r_p}(p\tau)$, then we must have $n_\infty=0$ by comparing the constant terms. Hence $r_1=-pr_p$ and $k=(1-p)r_p/2$. Now Corollary \ref{coro:noCuspForm} implies there are just two cases:
\begin{enumerate}
\item $p=2$ and $2 < k < 4(1+\frac{1}{m_1})$,
\item $p=3$ and $2 < k < 3(1+\frac{1}{m_1})$.
\end{enumerate}
\begin{thm}
\label{thm:EisEquEtaInfty}
If $n_1$ is any rational number with $1/2 < n_1 \leq 1$, then we have
\begin{equation*}
E_{4n_1,2}^{\rmi\infty}(\tau;16n_1,-8n_1)=\eta^{16n_1}(\tau)\eta^{-8n_1}(2\tau).
\end{equation*}
On the other hand, if $n_1$ is a rational number with $ 2/3 < n_1 \leq 1$, then we have
\begin{equation*}
E_{3n_1,3}^{\rmi\infty}(\tau;9n_1,-3n_1)=\eta^{9n_1}(\tau)\eta^{-3n_1}(3\tau).
\end{equation*}
\end{thm}
\begin{proof}
In both cases, the usage of the symbol $n_1$ as an arbitrarily chosen rational coincides with the usage $n_1=(pr_1+r_p)/24$ before. By Lemma \ref{lemm:etaQuoChar} and Lemma \ref{lemm:ordEtaQuo}, we have $\eta^{16n_1}(\tau)\eta^{-8n_1}(2\tau) \in M_{4n_1}(\widetilde{\Gamma_0(2)^{2m_1}}; \chi_{16n_1,-8n_1})$ and $\eta^{9n_1}(\tau)\eta^{-3n_1}(3\tau) \in M_{3n_1}(\widetilde{\Gamma_0(3)^{2m_1}}; \chi_{9n_1,-3n_1})$ whenever $n_1 \in \numQ_{>0}$. It follows from Theorem \ref{thm:mainEis} that, if $1/2 < n_1$ then $E_{4n_1,2}^{\rmi\infty}(\tau;16n_1,-8n_1)$ is well-defined and belongs to $M_{4n_1}(\widetilde{\Gamma_0(2)^{2m_1}}; \chi_{16n_1,-8n_1})$, and if $2/3 < n_1$ then $E_{3n_1,3}^{\rmi\infty}(\tau;9n_1,-3n_1)$ is well-defined and belongs to $M_{3n_1}(\widetilde{\Gamma_0(3)^{2m_1}}; \chi_{9n_1,-3n_1})$. Put $f=E_{4n_1,2}^{\rmi\infty}(\tau;16n_1,-8n_1)-\eta^{16n_1}(\tau)\eta^{-8n_1}(2\tau)$ and $g=E_{3n_1,3}^{\rmi\infty}(\tau;9n_1,-3n_1)-\eta^{9n_1}(\tau)\eta^{-3n_1}(3\tau)$. The constant terms of $f$ and $g$ vanish since the coefficient of the leading term of the $q$-expansion of any (rational-weight) eta-qoutient must be $1$ (Remark \ref{rema:etaCoeff}). Moreover, the constant terms of $f\vert_{4n_1}\widetilde{\gamma_1}$ and $g\vert_{3n_1}\widetilde{\gamma_1}$ are also zero due to Theorem \ref{thm:EisMod} and Lemma \ref{lemm:ordEtaQuo}. Hence $f$ under the condition $1/2 < n_1$ and $g$ under the condition $2/3 < n_1$ are both cusp forms. The desired identities follow from this observation and Corollary \ref{coro:noCuspForm}.
\end{proof}

Now we deal with $E_{k,p}^{1}(\tau;r_1,r_p)$ in a similar manner. Recall $\gamma_1=\tbtmat{1}{0}{1}{1}$ which satisfies $\gamma_1(\rmi\infty)=1/1$. By Theorem \ref{thm:EisMod} we know that $\lim_{\Im\tau\to+\infty}E_{k,p}^{1}(\tau;r_1,r_p)\vert_k\widetilde{\gamma_1}=1$. So we should search for $\eta^{r_1}(\tau)\eta^{r_p}(p\tau)$ with the property $\lim_{\Im\tau\to+\infty}\eta^{r_1}(\tau)\eta^{r_p}(p\tau)\vert_k\widetilde{\gamma_1}$ is nonzero. This means $n_1=0$ by Lemma \ref{lemm:ordEtaQuo}. Hence $r_p=-pr_1$ and $k=(1-p)r_1/2$. Now Corollary \ref{coro:noCuspForm} implies there are just two cases:
\begin{enumerate}
\item $p=2$ and $2 < k < 4(1+\frac{1}{m_\infty})$,
\item $p=3$ and $2 < k < 3(1+\frac{1}{m_\infty})$.
\end{enumerate}
\begin{thm}
\label{thm:EisEquEta1}
If $n_\infty$ is any rational number with $1/2 < n_\infty \leq 1$, then we have
\begin{equation*}
E_{4n_\infty,2}^{1}(\tau;-8n_\infty,16n_\infty)=2^{8n_\infty}\cdot\eta^{-8n_\infty}(\tau)\eta^{16n_\infty}(2\tau).
\end{equation*}
On the other hand, if $n_\infty$ is a rational number with $ 2/3 < n_\infty \leq 1$, then we have
\begin{equation*}
E_{3n_\infty,3}^{1}(\tau;-3n_\infty,9n_\infty)=3^{9n_\infty/2}\etp{-\frac{n_\infty}{4}}\cdot\eta^{-3n_\infty}(\tau)\eta^{9n_\infty}(3\tau).
\end{equation*}
\end{thm}
\begin{proof}
Omitted since it is similar to the proof of Theorem \ref{thm:EisEquEtaInfty}. However there is one thing worth mentioning: where are the factors $2^{8n_\infty}$ and $3^{9n_\infty/2}\etp{-\frac{n_\infty}{4}}$ come from? They come from comparing the constant terms of the $q$-series of, say in the case $p=3$, $E_{3n_\infty,3}^{1}(\tau;-3n_\infty,9n_\infty)\vert_k\widetilde{\gamma_1}$ and $\eta^{-3n_\infty}(\tau)\eta^{9n_\infty}(3\tau)\vert_k\widetilde{\gamma_1}$. We know from Theorem \ref{thm:EisMod} that the constant term of $E_{k,p}^{1}(\tau;r_1,r_p)\vert_k\widetilde{\gamma_1}$ is $1$. On the other hand, we have
\begin{align*}
\eta^{r_1}(\tau)\eta^{r_p}(p\tau)\vert_{(r_1+r_p)/2}\widetilde{\tbtmat{1}{0}{1}{1}} &=\eta^{r_1}(\tau)\vert_{r_1/2}\widetilde{\tbtmat{1}{0}{1}{1}}\cdot\eta^{r_p}(p\tau)\vert_{r_p/2}\widetilde{\tbtmat{1}{0}{1}{1}}\\
&=\left(\etp{-\frac{1}{24}r_1}\eta^{r_1}(\tau)\right)\cdot\left(p^{-r_p/2}\etp{\frac{(p-3)r_p}{24}}\eta^{r_p}\left(\frac{\tau+1}{p}\right)\right).
\end{align*}
Consequently the term with the least exponent of the $q$-expansion of $\eta^{r_1}(\tau)\eta^{r_p}(p\tau)\vert_{(r_1+r_p)/2}\widetilde{\gamma_1}$ is
\begin{equation*}
p^{-r_p/2}\etp{-\frac{r_1}{24}+\frac{(p^2-3p+1)r_p}{24p}}q^{n_1/p}
\end{equation*}
where $n_1=(pr_1+r_p)/24$. This is exactly $2^{-8n_\infty}$ if we set $p=2$, $r_1=-8n_\infty$, $r_2=16n_\infty$, and is $3^{-9n_\infty/2}\etp{\frac{n_\infty}{4}}$ if we set $p=3$, $r_1=-3n_\infty$, $r_3=9n_\infty$.
\end{proof}

We now turn to an application. Recall the Gamma function $\Gamma(z)$ is defined by
\begin{equation*}
\frac{1}{\Gamma(z)}=z\rme^{\gamma z}\prod_{\nu \in \numgeq{Z}{1}}\left(1+\frac{z}{\nu}\right)\rme^{-z/\nu}, \quad z \in \numC,
\end{equation*}
where $\gamma=\lim_{n\to+\infty}(\sum_{\nu=1}^{n}1/\nu-\log n)$ is Euler's constant. By comparing Fourier coefficients of both sides of identities in Theorem \ref{thm:EisEquEtaInfty} and \ref{thm:EisEquEta1}, we obtain certain series expression of $\Gamma(k)$ for rational number $k$ with $2<k\leq4$. Since $\Gamma(z+1)=z\Gamma(z)$ for $z \neq 0,-1,-2,\dots$, we actually obtain expressions for any rational $k \neq 0,-1,-2,\dots$. First introduce some notation. We write
\begin{equation*}
\eta^{r_1}(\tau)\eta^{r_p}(p\tau)=q^{n_\infty}\sum_{n \in \numgeq{Z}{0}}A_p(n;r_1,r_p)q^n.
\end{equation*}
We have explained in Remark \ref{rema:etaCoeff} that $A_p(0;r_1,r_p)=1$. Since the $D$-th power of $\eta^{r_1}(\tau)\eta^{r_p}(p\tau)$ is an eta-quotient of integral powers, we have $A_p(n;r_1,r_p) \in \numQ$ for any $n \in \numgeq{Z}{0}$, $r_1,r_p \in \numQ$.
\begin{thm}
\label{thm:GammaValue}
Let $k$ be a rational number with $2<k\leq4$, and $n \in \numgeq{Z}{1}$ with $A_2(n;4k,-2k) \neq 0$. We have
\begin{equation*}
\Gamma(k)=\frac{(2\uppi)^kn^{k-1}}{A_2(n;4k,-2k)}\sum_{c\in\numgeq{Z}{1}}\frac{1}{(2c)^k}\sum_{\twoscript{0 \leq d < 2c}{\gcd(2c,d)=1}}\etp{\frac{nd}{2c}-k\left(2s(-d,2c)-s(-d,c)\right)}.
\end{equation*}
On the other hand, if $2 < k \leq 3$, and $n \in \numgeq{Z}{1}$ satisfies $A_3(n;3k,-k) \neq 0$, then we have
\begin{equation*}
\Gamma(k)=\frac{(2\uppi)^kn^{k-1}}{A_3(n;3k,-k)}\sum_{c\in\numgeq{Z}{1}}\frac{1}{(3c)^k}\sum_{\twoscript{0 \leq d < 3c}{\gcd(3c,d)=1}}\etp{\frac{nd}{3c}-\frac{k}{2}\left(3s(-d,3c)-s(-d,c)\right)}.
\end{equation*}
\end{thm}
\begin{proof}
An immediate corollary of Theorem \ref{thm:EisEquEtaInfty}.
\end{proof}
Of course, there are another two series expansions of $\Gamma(k)$ that correspond to Theorem \ref{thm:EisEquEta1}. We omit them here since these series involve more complicated exponential sums of Dedekind sums.

It can be verified that $A_p(1;r_1,r_p)=-r_1$ for any rational $r_1, r_p$. Hence we have the following corollary.
\begin{coro}
\label{coro:GammaValue}
Let $k$ be a rational number with $2<k\leq4$. We have
\begin{equation}
\label{eq:GammaValuen1p2}
\Gamma(k)=-\frac{(2\uppi)^k}{4k}\sum_{c\in\numgeq{Z}{1}}\frac{1}{(2c)^k}\sum_{\twoscript{0 \leq d < 2c}{\gcd(2c,d)=1}}\etp{\frac{d}{2c}-k\left(2s(-d,2c)-s(-d,c)\right)}.
\end{equation}
On the other hand, if $2 < k \leq 3$, then we have
\begin{equation}
\label{eq:GammaValuen1p3}
\Gamma(k)=-\frac{(2\uppi)^k}{3k}\sum_{c\in\numgeq{Z}{1}}\frac{1}{(3c)^k}\sum_{\twoscript{0 \leq d < 3c}{\gcd(3c,d)=1}}\etp{\frac{d}{3c}-\frac{k}{2}\left(3s(-d,3c)-s(-d,c)\right)}.
\end{equation}
\end{coro}
\begin{rema}
By a continuity argument, actually the identity \eqref{eq:GammaValuen1p2} holds for any \emph{real} number $k$ with $2 < k \leq 4$, and \eqref{eq:GammaValuen1p3} holds for any \emph{real} number $k$ with $2 < k \leq 3$. It is interesting to study how these identities can be extended to imaginary values of $k$.
\end{rema}

\section{Discussion: special exponents}
\label{sec:Discussion}
In this section, we discuss how Eisenstein series $E_{k,p}^{\rmi\infty}(\tau;r_1,r_p)$, the central functions of this paper (see \eqref{eq:Eisinfty} or \eqref{eq:FourierExpEisInfty}), may be reduced to well-known functions if $r_1, r_p$ assume some special values. For this purpose, it is sometimes useful to rewrite \eqref{eq:FourierExpEisInfty} as
\begin{multline}
\label{eq:EisFourierInftyAnother}
E_{k,p}^{\rmi\infty}(\tau;r_1,r_p)=1+\etp{-\frac{k}{4}}\frac{(2\uppi)^k}{\Gamma(k)}\sum_{n \in \numgeq{Z}{1}}n^{k-1}\sum_{c \in \numgeq{Z}{1}}\frac{1}{(pc)^k}\times \\
\sum_{\twoscript{0 \leq d < pc}{\gcd(pc, d)=1}}\etp{\frac{nd}{pc}}\chi_{r_1,r_p}\widetilde{\tbtmat{a}{b}{pc}{d}}^{-1}q^n,
\end{multline}
where for each pair $(c,d)$, $a$ and $b$ are two integers satisfying $ad-pcb=1$. Recall that we require $\chi_{r_1,r_p}(\widetilde{T})=1$.

In the simplest case $\chi_{r_1,r_p}\widetilde{\tbtmat{a}{b}{pc}{d}}=1$ for any $c \in \numgeq{Z}{1}$ and $0 \leq d < pc$ with $\gcd(pc,d)=1$, the sum for $d$ is Ramanujan's sum. The second-simplest case is the one in which $\chi_{r_1,r_p}\widetilde{\tbtmat{a}{b}{pc}{d}}=\chi(d)$ with $\chi$ some Dirichlet character modulo $p$. In this case, the sum for $d$ is a Gauss sum associated with a Dirichlet character. In both cases, the Fourier coefficients can be simplified to certain elementary arithmetical functions. For eta-quotients of integral powers, Gordon, Hughes \cite{GH93} and Newman \cite{New57}, \cite{New59} have given a criterion on when their multiplier systems are induced by Dirichlet characters. See also \cite[Proposition 5.9.2]{CoS17} for a proof. Here we will give a criterion on when $\chi_{r_1,r_p}$ is ``trivial'' for rational $r_1$ and $r_p$. For this purpose, we need some formulas, which do not involve Dedekind sums, for evaluating certain values $\chi_{r_1,r_p}(\widetilde{\gamma})$.
\begin{lemm}
\label{lemm:specialValuesChi}
Let $p$ be a prime and $r_1, r_p$ be two rationals. Let $D$ be a positive integer such that $D\cdot r_1$ and $D\cdot r_p$ are integers. Let $\chi_{r_1, r_p}\colon \widetilde{\Gamma_0(p)^{2D}}\to\numC^\times$ be the multiplier system of $\eta^{r_1}(\tau)\eta^{r_p}(p\tau)$ (which is a group character by the choice of $D$). The following formulas hold:
\begin{enumerate}
\item If $t \in \numgeq{Z}{1}$ is a factor of $p-1$, then
\begin{equation*}
\chi_{r_1, r_p}\widetilde{\tbtMat{-t}{-1}{p}{(p-1)/t}}=\etp{(r_1+r_p)\cdot\frac{p-(t^2+3t+1)}{24t}}.
\end{equation*}
\item If $t \in \numgeq{Z}{1}$ is a factor of $p+1$, then
\begin{equation*}
\chi_{r_1, r_p}\widetilde{\tbtMat{t}{1}{p}{(p+1)/t}}=\etp{(-r_1+r_p)\cdot\frac{p+t^2-3t+1}{24t}}.
\end{equation*}
\item If $t \in \numgeq{Z}{1}$ is a factor of $p-2$ and $p \geq 3$, then
\begin{equation*}
\chi_{r_1, r_p}\widetilde{\tbtMat{-t(p+1)/2}{-(p-1)/2}{p}{(p-2)/t}}=\etp{r_1\cdot\frac{p-(3t^2+6t+2)}{48t}+r_p\cdot\frac{(2-t^2)p-(t^2+6t+4)}{48t}}.
\end{equation*}
\item If $t \in \numgeq{Z}{1}$ is a factor of $p+2$ and $p \geq 3$, then
\begin{equation*}
\chi_{r_1, r_p}\widetilde{\tbtMat{t(p+1)/2}{(p+3)/2}{p}{(p+2)/t}}=\etp{r_1\cdot\frac{-p-(t^2-6t+2)}{48t}+r_p\cdot\frac{(t^2+2)p+t^2-6t+4}{48t}}.
\end{equation*}
\end{enumerate}
\end{lemm}
\begin{proof}
These formulas can be verified by a tedious calculation using Lemma \ref{lemm:etaQuoChar} and the following formulas for Dedekind sums:
\begin{align*}
s\left(\frac{c-1}{t},c\right)&=-\frac{(c-1)(c-(t^2+1))}{12tc},&\quad &c,t \in \numgeq{Z}{1},\, t\mid c-1; \\
s\left(\frac{c+1}{t},c\right)&=\frac{c^2+(t^2-6t+2)c+t^2+1}{12tc},&\quad &c,t \in \numgeq{Z}{1},\, t\mid c+1; \\
s\left(\frac{c-2}{t},c\right)&=-\frac{c^2-(2t^2+4)c+t^2+4}{24tc},&\quad &c,t \in \numgeq{Z}{1},\, t\mid c-2,\, 2 \nmid c; \\
s\left(\frac{c+2}{t},c\right)&=\frac{c^2+(2t^2-12t+4)c+t^2+4}{24tc},&\quad &c,t \in \numgeq{Z}{1},\, t\mid c+2,\, 2 \nmid c.
\end{align*}
These formulas for Dedekind sums themselves can be verified using the reciprocity law \cite[Theorem 3]{RW41}, or can be deduced from \cite[Chapter 3, Exercise 14,15,16,17]{Apo90}.
\end{proof}
\begin{prop}
\label{prop:whenChiTrivial}
Notations as in Lemma \ref{lemm:specialValuesChi}. Suppose that $\chi_{r_1,r_p}(\widetilde{T})=1$, that is to say, $r_1+pr_p \in 24\numZ$. Consider the following three statements:
\begin{enumerate}
\item For any $\tbtmat{a}{b}{c}{d} \in \Gamma_0(p)$ with $c > 0$ and $0 \leq d < c$, we have $\chi_{r_1,r_p}\widetilde{\tbtmat{a}{b}{c}{d}}=1$.
\item For any $\gamma \in \widetilde{\Gamma_0(p)^{2D}}$ we have $\chi_{r_1,r_p}(\gamma)=1$.
\item $r_1 \in 2\numZ$, $r_1+r_p \in 4\numZ$ and $pr_1+r_p \in 24\numZ$.
\end{enumerate}
We have $(1)\iff(2)\impliedby(3)$. Moreover, if $3 \mid p+1$ or $4 \mid p+1$ or $5 \mid p+1$, then $(2)\implies(3)$.
\end{prop}
\begin{proof}
For simplicity, let $\chi$ denote $\chi_{r_1,r_p}$.

$(1)\implies(2)$. Since $\widetilde{-I}=\widetilde{\tbtmat{1}{0}{p}{1}}\widetilde{\tbtmat{-1}{-1}{p}{p-1}}\widetilde{\tbtmat{1}{-1}{0}{1}}$, we have $\chi(\widetilde{-I})=1$. Let $\gamma=\tbtmat{a}{b}{c}{d} \in \Gamma_0(p)$ be arbitrary. If $c=0$, it follows from the fact $\chi_{r_1,r_p}(\widetilde{T})=1$ and $\chi(\widetilde{-I})=1$ that $\chi(\widetilde{\gamma})=1$. If $c>0$, there exist $q \in \numZ$ and $0 \leq r < c$ such that $d=cq+r$. Thus
\begin{equation*}
\chi\widetilde{\tbtMat{a}{b}{c}{d}}=\chi\widetilde{\tbtMat{a}{b-aq}{c}{d-cq}}\cdot\chi\widetilde{\tbtMat{1}{q}{0}{1}}=1.
\end{equation*}
If $c < 0$, we have $\chi\widetilde{\tbtmat{a}{b}{c}{d}}=\chi\widetilde{\tbtmat{-a}{-b}{-c}{-d}}=1$ since $\widetilde{\tbtmat{a}{b}{c}{d}}\cdot\widetilde{(-I)}=\widetilde{\tbtmat{-a}{-b}{-c}{-d}}$ and $\chi(\widetilde{-I})=1$. Now let $\gamma=\elesltRaD{a}{b}{c}{d}{\etp{\frac{j}{2D}}}{2D} \in \widetilde{\Gamma_0(p)^{2D}}$ be arbitrary, where $\tbtmat{a}{b}{c}{d}\in\Gamma_0(p)$ and $0 \leq j < 2D$. We have
\begin{equation*}
\chi(\gamma)=\chi\widetilde{\tbtMat{a}{b}{c}{d}}\cdot\chi\left(\widetilde{-I}^{2j}\right)=1.
\end{equation*}
Hence $(1)\implies(2)$.

$(2)\implies(1)$. This is obvious.

$(3)\implies(1)$. We need a theorem of Ligozat \cite{Lig75}:
\begin{equation}
\label{eq:ligozat}
\eta\vert_{1/2}\widetilde{\tbtmat{a}{b}{c}{d}}=\legendre{c}{a}\cdot\etp{\frac{1}{24}(a(b-c+3)-3)}\cdot\eta,\quad \tbtmat{a}{b}{c}{d}\in\slZ,\, c\geq0, \gcd(a,6)=1,
\end{equation}
where $\legendre{c}{a}$ is the Kronecker-Jacobi symbol. If $r_1,r_p \in \numZ$, $r_1+r_p \in 4\numZ$, $r_1+pr_p \in 24\numZ$ and $pr_1+r_p \in 24\numZ$, we can deduce from Ligozat's theorem that, if $p \mid c$ then
\begin{equation}
\label{eq:gcda6}
\eta^{r_1}(\tau)\eta^{r_p}(p\tau)\vert_{k'}\widetilde{\tbtmat{a}{b}{c}{d}}=\legendre{P}{a}\cdot\eta^{r_1}(\tau)\eta^{r_p}(p\tau),
\end{equation}
where $k'=(r_1+r_p)/2$, $P=p$ if $2 \nmid r_1$ and $P=1$ if $2 \mid r_1$. Since $2 \mid r_1$ by assumption, we have $\chi\widetilde{\tbtmat{a}{b}{c}{d}}=1$ for $\tbtmat{a}{b}{c}{d}\in\Gamma_0(p)$ with $c \geq 0$ and $\gcd(a,6)=1$. Since such matrices generate $\Gamma_0(p)$, the statement $(1)$ follows.

$(2)\implies(3)$. It follows from the definition that $\chi(\widetilde{-I})=\etp{-k'/2}$. However, $\chi(\widetilde{-I})=1$, and so $k'\in2\numZ$, that is, $r_1+r_p \in 4\numZ$. We have known from the proof of Corollary \ref{coro:noCuspForm} that $\chi(\widetilde{\gamma_1}\widetilde{T}^p\widetilde{\gamma_1}^{-1})=\etp{n_1}$, where $n_1=(pr_1+r_p)/24$. Hence $pr_1+r_p \in 24\numZ$. It remains to prove $r_1 \in 2\numZ$, which is the most difficult part of this proof and is the only part relying on the assumption $3 \mid p+1$ or $4 \mid p+1$ or $5 \mid p+1$. First assume that $3 \mid p+1$. By Lemma \ref{lemm:specialValuesChi}(2) we have
\begin{align*}
\chi\widetilde{\tbtMat{3}{1}{p}{(p+1)/3}}&=\etp{(-r_1+r_p)\cdot\frac{p+1}{72}}, \\
\chi\widetilde{\tbtMat{1}{1}{p}{p+1}}&=\etp{(-r_1+r_p)\cdot\frac{p-1}{24}}.
\end{align*}
It follows that $(-r_1+r_p)(p+1) \in 72\numZ$ and $(-r_1+r_p)(p-1) \in 24\numZ$. Cancelling $p$ gives that $-r_1+r_p \in 12\numZ$. This and $r_1+r_p \in 4\numZ$ together imply $r_1 \in 2\numZ$. The proof for the case $4 \mid p+1$ is similar. In the last case, $5 \mid p+1$, this method only leads to $r_1 \in \numZ$. Since $r_1$ and $r_p$ are both integers now, we can apply formula \eqref{eq:gcda6} again and obtain $r_1 \in 2\numZ$.
\end{proof}
\begin{conj}
For any prime $p$, the three statements in the Proposition \ref{prop:whenChiTrivial} are equivalent.
\end{conj}

If the condition (1) in Proposition \ref{prop:whenChiTrivial} holds, then the Eisenstein series $E_{k,p}^{\rmi\infty}(\tau;r_1,r_p)$ can be simplified to
\begin{equation*}
E_{k,p}^{\rmi\infty}(\tau;r_1,r_p)=1-\frac{2k}{B_k}\frac{1}{p^k-1}\sum_{n \in \numgeq{Z}{1}}\left(p^k\sigma_{k-1}\left(n/p\right)-\sigma_{k-1}(n)\right)q^n,
\end{equation*}
where $B_k$ is Bernoulli numbers defined by $x/(\rme^x-1)=\sum_{k\geq 0}B_kx^k/k!$, and $\sigma_{k-1}(n)=\sum_{0 < d \mid n}d^{k-1}$ if $n \in \numgeq{Z}{1}$ and $\sigma_{k-1}(n)=0$ otherwise. This is the well-known Eisenstein series of even weight $k=4,6,8,\dots$, on the group $\Gamma_0(p)$ with trivial multiplier system. So Proposition \ref{prop:whenChiTrivial} tells us that such simple situation can not happen, at least for the case $3 \mid p+1$ or $4 \mid p+1$ or $5 \mid p+1$, if $r_1$ or $r_p$ are indeed fractions.

Next we consider another sort of special values of $r_1$ and $r_p$. If $r_1$ and $r_p$ are both integers, then there are several formulas for $\chi_{r_1, r_p}$ that does not involve Dedekind sums. One of these formulas is based on Petersson's transformation equation of $\eta$. See \cite{Kno70} or \cite[Eq.(15)]{ZZ21}. Since we have restricted ourselves to integer exponents, we may choose $D=1$, that is, we may use double covers of modular groups. We put $\widetilde{\Gamma_0(p)}=\widetilde{\Gamma_0(p)^2}$ for simplicity.
\begin{lemm}
\label{lemm:r1rpIntegers}
Let $p$ be a prime and $r_1,r_p$ be integers with $n_\infty=(r_1+pr_p)/24 \in \numZ$. Let $\tbtmat{a}{b}{c}{d} \in \Gamma_0(p)$ be arbitrary. If $p \geq 5$, then
\begin{equation*}
\chi_{r_1,r_p}\widetilde{\tbtMat{a}{b}{c}{d}}=\legendre{d}{p}^{r_p}.
\end{equation*}
If $p=3$, then
\begin{equation*}
\chi_{r_1,r_3}\widetilde{\tbtMat{a}{b}{c}{d}}=\legendre{d}{3}^{r_3}\etp{-r_3\cdot\frac{c(a+d)}{9}}.
\end{equation*}
If $p=2$, then
\begin{equation*}
\chi_{r_1,r_2}\widetilde{\tbtMat{a}{b}{c}{d}}=\legendre{c/2}{a}^{r_2}\etp{r_2\cdot\frac{a(c/2-1)+1}{8}}.
\end{equation*}
\end{lemm}
\begin{proof}
The formula for the case $p \geq 5$ can be deduced immediately from a result of Gordon, Hughes and Newman (see also \cite[Proposition 5.9.2]{CoS17}). But we present another proof here since it is also valid for the case $p=3$. Put $n_1=(pr_1+r_p)/24$. We shall prove that, if $p \geq 3$ then
\begin{equation}
\label{eq:chir1rpIntegers}
\chi_{r_1,r_p}\widetilde{\tbtMat{a}{b}{c}{d}}=\legendre{d}{p}^{r_p}\etp{\frac{n_1(a+d-bdc-3)c}{p}},\quad\tbtMat{a}{b}{c}{d}\in\Gamma_0(p).
\end{equation}
We prove by cases. In the case $2 \nmid c$, we have
\begin{equation*}
\chi_{r_1,r_p}\widetilde{\tbtMat{a}{b}{c}{d}}=\legendre{d}{\abs{c}}^{r_1+r_p}\legendre{d}{p}^{r_p}\etp{\frac{n_1(a+d-bdc-3)c}{p}}
\end{equation*}
by \cite[Eq.(15)]{ZZ21} and the fact $n_\infty \in \numZ$. Again from $n_\infty \in \numZ$ and $p \geq 3$ we know that $2\mid r_1+r_p$, and so \eqref{eq:chir1rpIntegers} follows in the case $2 \nmid c$. On the other hand, if $2 \mid c$, then
\begin{equation}
\label{eq:eq:chir1rpIntegersCasecEven}
\chi_{r_1,r_p}\widetilde{\tbtMat{a}{b}{c}{d}}=\legendre{c}{d}^{r_1+r_p}\legendre{p}{d}^{r_p}\etp{\frac{n_1(a-2d-bdc)c}{p}}\etp{\frac{r_1+r_p}{4}\frac{d-1}{2}}
\end{equation}
still by \cite[Eq.(15)]{ZZ21} and $n_\infty \in \numZ$. Since $2\mid r_1+r_p$, we have $\legendre{c}{d}^{r_1+r_p}=1$. It follows from the fact $2 \mid p/c$ and $2 \mid r_1+r_p$ that $\etp{\frac{n_1(a-2d)c}{p}}=\etp{\frac{n_1(a+d-3)c}{p}}$. By the Jacobi reciprocity law\footnote{Many authors state Jacobi reciprocity law for positive odd integers. It is useful to admit negative odd integers. (So we should use Jacobi-Kronecker symbol.) The formula reads $\legendre{Q}{P}\cdot\legendre{P}{Q}=\varepsilon(P,Q)\cdot(-1)^{(P-1)(Q-1)/4}$, where $P$ and $Q$ are any coprime odd integers, and $\varepsilon(P,Q)=-1$ if both $P$, $Q$ are negative, $\varepsilon(P,Q)=1$ otherwise.} we deduce that
\begin{equation*}
\legendre{p}{d}^{r_p}\cdot\etp{\frac{r_1+r_p}{4}\frac{d-1}{2}} = \legendre{d}{p}^{r_p}\cdot\etp{\frac{r_1+pr_p}{4}\frac{d-1}{2}}=\legendre{d}{p}^{r_p}.
\end{equation*}
Inserting these identities into \eqref{eq:eq:chir1rpIntegersCasecEven} gives \eqref{eq:chir1rpIntegers} in the case $2 \mid c$.
Now if $p \geq 5$, then $n_\infty \in \numZ$ implies that $n_1 \in \numZ$ since $n_1=pn_\infty-r_p(p^2-1)/24$. Thus, it follows from \eqref{eq:chir1rpIntegers} that $\chi_{r_1,r_p}\widetilde{\tbtmat{a}{b}{c}{d}}=\legendre{d}{p}^{r_p}$. The assertion in the case $p=3$ follows similarly. Finally, we aim to prove the case $p=2$. It follows from \eqref{eq:ligozat} that
\begin{equation*}
\chi_{r_1,r_p}\widetilde{\tbtMat{a}{b}{c}{d}}=\legendre{c}{a}^{r_1+r_p}\legendre{p}{a}^{r_p}\etp{\frac{(a-1)(r_1+r_p)}{8}}\etp{n_\infty ab-n_1ac/p}
\end{equation*}
where $\tbtmat{a}{b}{c}{d}\in\Gamma_0(p)$ with $c > 0$ and $\gcd(a,6)=1$. Setting $p=2$ and using the facts $n_1=pn_\infty-r_p(p^2-1)/24$, $n_\infty \in \numZ$, we find that the desired assertion holds for $\tbtmat{a}{b}{c}{d}\in\Gamma_0(2)$ with $c > 0$ and $\gcd(a,6)=1$. If $\gcd(a,6) \neq 1$, then $\gcd(a+2c,6)=1$, and
\begin{equation*}
\legendre{c/2}{a+2c}=\legendre{2c}{a+2c}=\legendre{2c}{a}=\legendre{c/2}{a}
\end{equation*}
since $a$ is odd. Note that $\chi_{r_1,r_p}(\widetilde{T})=1$. Therefore
\begin{equation*}
\chi_{r_1,r_p}\widetilde{\tbtMat{a}{b}{c}{d}}=\chi_{r_1,r_p}\widetilde{\tbtMat{a+2c}{b+2d}{c}{d}}=\legendre{c/2}{a}^{r_2}\etp{r_2\cdot\frac{a(c/2-1)+1}{8}}.
\end{equation*}
This means the desired assertion holds for any $\tbtmat{a}{b}{c}{d}\in\Gamma_0(p)$ with $c > 0$, whatever $\gcd(a,6)$ is. It also holds for $c=0$, as one can directly verify. If $c<0$, then
\begin{equation*}
\chi_{r_1,r_p}\widetilde{\tbtMat{a}{b}{c}{d}}=\chi_{r_1,r_p}\widetilde{\tbtMat{-a}{-b}{-c}{-d}}\cdot\chi_{r_1,r_p}(\widetilde{-I})^{-1}=\legendre{-c/2}{-a}^{r_2}\etp{r_2\cdot\frac{ac/2+a-1}{8}}.
\end{equation*}
The desired assertion follows from this and the fact $\legendre{-c/2}{-a}^{r_2}\etp{r_2\cdot\frac{a-1}{4}}=\legendre{c}{a}^{r_2}$. This concludes the proof for the case $p=2$, hence the whole proof.
\end{proof}

Look at the expression \eqref{eq:EisFourierInftyAnother}. If $r_1$ and $r_p$ are both integers, then we can plug Lemma \ref{lemm:r1rpIntegers} into this expression and simplify the Fourier coefficients. For $p \geq 5$, this leads to well-known Eisenstein series of integral weights (greater than $2$) on the groups $\Gamma_0(p)$ with real Dirichlet characters. The cases $p=2$ and $p=3$ are more interesting. For $p=3$, the Fourier coefficients involve Kloosterman sums with characters, which are defined by
\begin{equation*}
K(\psi,m,n;c)=\sum_{\twoscript{r \bmod c}{\gcd(r,c)=1}}\psi(r)\etp{\frac{mr+nr^{-1}}{c}}.
\end{equation*}
\begin{prop}
\label{prop:Eisp3}
If $r_1$ and $r_3$ are integers with $r_1+3r_3$ divisible by $24$. Let $k$ be an integer greater than $2$ satisfying $k \equiv r_3 \bmod 2$. Then we have
\begin{equation*}
E_{k,3}^{\rmi\infty}(\tau;r_1,r_3)=1+\etp{-\frac{k}{4}}\frac{(2\uppi)^k}{(k-1)!}\sum_{n \in \numgeq{Z}{1}}n^{k-1}\sum_{c \in \numgeq{Z}{1}}\frac{K(\legendre{\bullet}{3}^{r_3},c^2r_3+n,c^2r_3;3c)}{(3c)^k}q^n,
\end{equation*}
where $\legendre{\bullet}{3}^{r_3}$ means the function sending $d$ to $\legendre{d}{3}^{r_3}$. These functions depend only on $r_3 \bmod 6$, and $k$, so we may let $r_3 \in \{0,-1,-2,-3,-4,-5\}$, $r_1=-3r_3$ and lose nothing.
\end{prop}
\begin{proof}
The identity is a combination of \eqref{eq:EisFourierInftyAnother}, Lemma \ref{lemm:r1rpIntegers} and the definition of Kloosterman sums. The weight $k$ is required to be in $(r_1+r_3)/2+2\numZ$ (the paragraph before Theorem \ref{thm:mainEis}). Since $24 \mid r_1+3r_3$, we have $(r_1+r_3)/2 \equiv r_3 \bmod 2$. Thus, the requirement of $k$ is equivalent to that $k \in r_3+2\numZ$. From the formula for $E_{k,3}^{\rmi\infty}(\tau;r_1,r_3)$ just verified, we see that the functions are independent of $r_1$, and depend on $r_3$ only modulo $6$.
\end{proof}
\begin{examp}
If $3 \mid r_3$, then $E_{k,3}^{\rmi\infty}(\tau;r_1,r_3)$ becomes Eisenstein series whose multiplier systems are induced by real Dirichlet characters. For example, $r_3=-3$ is in this case. This case is also included in Theorem \ref{thm:EisEquEtaInfty}. So we have $E_{3,3}^{\rmi\infty}(\tau;9,-3)=\eta^9(\tau)\eta^{-3}(3\tau)$. After a standard simplification of \eqref{eq:EisFourierInftyAnother}, we obtain
\begin{equation*}
\eta^9(\tau)\eta^{-3}(3\tau)=1-9\sum_{n \in \numgeq{Z}{1}}\sum_{0 < m \mid n}\legendre{m}{3}m^2q^n.
\end{equation*}
This identity was first proved by Carlitz \cite[Eq.(3.1)]{Car53}. See also \cite[Lemma 14.2.5]{AB05}.
\end{examp}

It is not difficult to prove that, some functions in Proposition \ref{prop:Eisp3}, namely, all $E_{k,3}^{\rmi\infty}(\tau;3k,-k)$ with $k=3,4,5,\dots$, can be expressed as a linear combination of the eta-quotients $\eta^{3k}(\tau)\eta^{-k}(3\tau)$, $\eta^{3k-12}(\tau)\eta^{-k+12}(3\tau)$, $\eta^{3k-24}(\tau)\eta^{-k+24}(3\tau),\dots$ using \cite[lemma 6.4]{ZZ21}. In fact $k$ is allowed to assume rational values, which may be seen as a generalization of the second identity in Theorem \ref{thm:EisEquEtaInfty}. It is an interesting problem to determine the coefficients explicitly.

We omit discussion on the case $p=2$.

%

\bibliographystyle{elsarticle-num}
\bibliography{main}

\begin{thebibliography}{10}
\expandafter\ifx\csname url\endcsname\relax
  \def\url#1{\texttt{#1}}\fi
\expandafter\ifx\csname urlprefix\endcsname\relax\def\urlprefix{URL }\fi
\expandafter\ifx\csname href\endcsname\relax
  \def\href#1#2{#2} \def\path#1{#1}\fi

\bibitem{Sage}
SageMath, \href{https://www.sagemath.org}{The sage mathematics software system
  (version 9.3)}, The Sage Developers (2021).
\newline\urlprefix\url{https://www.sagemath.org}

\bibitem{Shi73}
G.~Shimura, \href{https://doi.org/10.2307/1970831}{On modular forms of half
  integral weight}, Ann. of Math. (2) 97 (1973) 440--481.
\newblock \href {https://doi.org/10.2307/1970831} {\path{doi:10.2307/1970831}}.
\newline\urlprefix\url{https://doi.org/10.2307/1970831}

\bibitem{Ibu00}
T.~Ibukiyama, \href{https://doi.org/10.1007/BF02940923}{Modular forms of
  rational weights and modular varieties}, Abh. Math. Sem. Univ. Hamburg 70
  (2000) 315--339.
\newblock \href {https://doi.org/10.1007/BF02940923}
  {\path{doi:10.1007/BF02940923}}.
\newline\urlprefix\url{https://doi.org/10.1007/BF02940923}

\bibitem{Ibu20_2}
T.~Ibukiyama, \href{https://doi.org/10.1007/s40993-019-0183-9}{Graded rings of
  modular forms of rational weights}, Res. Number Theory 6~(1) (2020) Paper No.
  8, 13.
\newblock \href {https://doi.org/10.1007/s40993-019-0183-9}
  {\path{doi:10.1007/s40993-019-0183-9}}.
\newline\urlprefix\url{https://doi.org/10.1007/s40993-019-0183-9}

\bibitem{Apo90}
T.~M. Apostol, \href{https://doi.org/10.1007/978-1-4612-0999-7}{Modular
  functions and {D}irichlet series in number theory}, 2nd Edition, Vol.~41 of
  Graduate Texts in Mathematics, Springer-Verlag, New York, 1990.
\newline\urlprefix\url{https://doi.org/10.1007/978-1-4612-0999-7}

\bibitem{Ono04}
K.~Ono, The web of modularity: arithmetic of the coefficients of modular forms
  and {$q$}-series, Vol. 102 of CBMS Regional Conference Series in Mathematics,
  Published for the Conference Board of the Mathematical Sciences, Washington,
  DC; by the American Mathematical Society, Providence, RI, 2004.

\bibitem{Rem98}
R.~Remmert, \href{https://doi.org/10.1007/978-1-4757-2956-6}{Classical topics
  in complex function theory}, Vol. 172 of Graduate Texts in Mathematics,
  Springer-Verlag, New York, 1998, translated from the German by Leslie Kay.
\newblock \href {https://doi.org/10.1007/978-1-4757-2956-6}
  {\path{doi:10.1007/978-1-4757-2956-6}}.
\newline\urlprefix\url{https://doi.org/10.1007/978-1-4757-2956-6}

\bibitem{ZZ21}
H.-G. Zhou, X.-J. Zhu, \href{https://doi.org/10.1016/j.jnt.2023.02.017}{Double
  coset operators and eta-quotients}, J. Number Theory (2023).
\newblock \href {http://arxiv.org/abs/2110.06768v1}
  {\path{arXiv:2110.06768v1}}.
\newline\urlprefix\url{https://doi.org/10.1016/j.jnt.2023.02.017}

\bibitem{RW41}
H.~Rademacher, A.~Whiteman, \href{https://doi.org/10.2307/2371532}{Theorems on
  {D}edekind sums}, Amer. J. Math. 63 (1941) 377--407.
\newblock \href {https://doi.org/10.2307/2371532} {\path{doi:10.2307/2371532}}.
\newline\urlprefix\url{https://doi.org/10.2307/2371532}

\bibitem{GH93}
B.~Gordon, K.~Hughes,
  \href{https://doi.org/10.1090/conm/143/01008}{Multiplicative properties of
  {$\eta$}-products. {II}}, in: A tribute to {E}mil {G}rosswald: number theory
  and related analysis, Vol. 143 of Contemp. Math., Amer. Math. Soc.,
  Providence, RI, 1993, pp. 415--430.
\newblock \href {https://doi.org/10.1090/conm/143/01008}
  {\path{doi:10.1090/conm/143/01008}}.
\newline\urlprefix\url{https://doi.org/10.1090/conm/143/01008}

\bibitem{New57}
M.~Newman, \href{https://doi.org/10.1112/plms/s3-7.1.334}{Construction and
  application of a class of modular functions}, Proc. London Math. Soc. (3) 7
  (1957) 334--350.
\newblock \href {https://doi.org/10.1112/plms/s3-7.1.334}
  {\path{doi:10.1112/plms/s3-7.1.334}}.
\newline\urlprefix\url{https://doi.org/10.1112/plms/s3-7.1.334}

\bibitem{New59}
M.~Newman, \href{https://doi.org/10.1112/plms/s3-9.3.373}{Construction and
  application of a class of modular functions. {II}}, Proc. London Math. Soc.
  (3) 9 (1959) 373--387.
\newblock \href {https://doi.org/10.1112/plms/s3-9.3.373}
  {\path{doi:10.1112/plms/s3-9.3.373}}.
\newline\urlprefix\url{https://doi.org/10.1112/plms/s3-9.3.373}

\bibitem{CoS17}
H.~Cohen, F.~Str\"omberg, Modular forms, Vol. 179 of Graduate Studies in
  Mathematics, American Mathematical Society, Providence, RI, 2017, a classical
  approach.

\bibitem{Lig75}
G.~Ligozat, Courbes modulaires de genre {$1$}, Suppl\'{e}ment au Bull. Soc.
  Math. France, Tome 103, no. 3, Soci\'{e}t\'{e} Math\'{e}matique de France,
  Paris, 1975, bull. Soc. Math. France, M\'{e}m. 43.

\bibitem{Kno70}
M.~I. Knopp, Modular functions in analytic number theory, Markham Publishing
  Co., Chicago, Ill., 1970.

\bibitem{Car53}
L.~Carlitz, \href{https://doi.org/10.1093/qmath/4.1.168}{Note on some partition
  formulae}, Quart. J. Math. Oxford Ser. (2) 4 (1953) 168--172.
\newblock \href {https://doi.org/10.1093/qmath/4.1.168}
  {\path{doi:10.1093/qmath/4.1.168}}.
\newline\urlprefix\url{https://doi.org/10.1093/qmath/4.1.168}

\bibitem{AB05}
G.~E. Andrews, B.~C. Berndt, Ramanujan's lost notebook. {P}art {I}, Springer,
  New York, 2005.

\end{thebibliography}
\end{document}